\documentclass[12pt]{amsart}

\usepackage{graphicx}

\usepackage{amsthm}
\usepackage{amsmath}
\usepackage{amssymb}
\usepackage{enumerate}

\usepackage{amsthm}
\usepackage{amsmath}
\usepackage{amssymb}
\usepackage{graphicx}
\usepackage{enumerate}
\usepackage{color}
\allowdisplaybreaks
\numberwithin{equation}{section}
\numberwithin{figure}{section}

\theoremstyle{plain}
\newtheorem{theorem}{\sffamily Theorem}[section]
\newtheorem{proposition}[theorem]{\sffamily Proposition}
\newtheorem{lemma}[theorem]{\sffamily Lemma}
\newtheorem{corollary}[theorem]{\sffamily Corollary}
\newtheorem{example}[theorem]{\sffamily Example}
\newtheorem{remark}[theorem]{\sffamily Remark}
\newtheorem{definition}[theorem]{\sffamily Definition}
\newtheorem{conjecture}{\sffamily Conjecture}

\def\BET{\begin{theorem}}
\def\ENT{\end{theorem}}
\def\BEP{\begin{proposition}}
\def\ENP{\end{proposition}}
\def\BEL{\begin{lemma}}
\def\ENL{\end{lemma}}
\def\BEC{\begin{corollary}}
\def\ENC{\end{corollary}}
\def\BEE{\begin{example}\rm}
\def\ENE{\end{example}}
\def\BER{\begin{remark} \rm}
\def\ENR{\end{remark}}
\def\BED{\begin{definition} \rm}
\def\END{\end{definition}}
\def\BECJ{\begin{conjecture}}
\def\ENCJ{\end{conjecture}}

\def\bea{\begin{eqnarray}}
\def\eea{\end{eqnarray}}

\def\beas{\begin{eqnarray*}}
\def\eeas{\end{eqnarray*}}

\def\beq{\begin{equation}}
\def\eeq{\end{equation}}

\def\beal{\begin{align*}}

\def\eeal{ \end{align*} }

\newcommand{\row}{ \nonumber \\ & & }
\newcommand{\roweq}{\nonumber \\ &=& }
\newcommand{\rowleq}{\nonumber \\  & \leq & }

\newcommand{\rowpl}{\nonumber \\ & \ \ &+  }
\newcommand{\rowmi}{\nonumber \\ & \ \ &-  }

\newcommand{\bfU}{{\bf U}}
\newcommand{\bfV}{{\bf V}}
\newcommand{\bfW}{{\bf W}}

\newcommand{\bfw}{{\bf w}}

\newcommand{\bbB}{{\mathbb B}}
\newcommand{\bbC}{{\mathbb C}}

\newcommand{\bbK}{{\mathbb K}}

\newcommand{\bbN}{{\mathbb N}}

\newcommand{\bbP}{{\mathbb P}}

\newcommand{\bbR}{{\mathbb R}}

\newcommand{\bbZ}{{\mathbb Z}}

\newcommand{\cH}{{\mathcal H}}

\newcommand{\cJ}{{\mathcal J}}
\newcommand{\cK}{{\mathcal K}}

\newcommand{\cP}{{\mathcal P}}

\newcommand{\cT}{{\mathcal T}}
\newcommand{\cU}{{\mathcal U}}
\newcommand{\cV}{{\mathcal V}}
\newcommand{\cW}{{\mathcal W}}
\newcommand{\cX}{{\mathcal X}}
\newcommand{\cY}{{\mathcal Y}}

\begin{document}

\title[water-wave problem in a shallow canal]{Band-gap structure of the spectrum of
the water-wave problem in a shallow canal with a periodic family of
deep pools}

\author{Sergei A. Nazarov}
\address{St. Petersburg State University,
Universitetskaya nab. 7--9, St. Petersburg, 199034, 
Russia, and
Institute for Problems in Mechanical Engineering of RAS, St. Petersburg, 199178, Russia 
}
\email{srgnazarov@yahoo.co.uk}

\author{Jari Taskinen}
\address{Department of Mathematics, University of Helsinki, 00014
Helsinki, Finland}
\email{jari.taskinen@helsinki.fi}


\thanks{The first named author was partially supported by  the Russian Foundation on Basic Research, project 18-01-00325.}



\begin{abstract}
We consider the linear water-wave problem in a periodic channel $\Pi^h \subset 
\bbR^3$, which is shallow except for a periodic array of deep potholes in it. 
Motivated by applications to 
surface wave propagation phenomena, we study the band-gap structure of the 
essential spectrum in the linear water-wave system, which includes 
the spectral Steklov boundary condition posed on the 
free water surface. We apply methods of asymptotic analysis, where the
most involved step is the construction and analysis of an appropriate boundary 
layer in a neighborhood of the joint of the potholes with the thin part
of the channel. Consequently, the existence of a spectral gap
for small enough $h$ is proven. 
\end{abstract}

\maketitle


\section{Introduction.}
\label{sec1}

\subsection{Formulation of the water-wave problem.} \label{sec1.1}
Let $x = (y,z) \in \bbR^2 \times \bbR$ be the Cartesian coordinate system
and let  $H:  \Omega^0 = \bbR \times (- \ell,\ell) \to \bbR$ 
be a smooth function, which is periodic in the variable $x_1 = y_1$ for
$y = (y_1,y_2) \in \Omega^0$ as well as positive; precisely, $H(y) \geq H_0 > 
0$ in $\overline{\Omega^0}$. By rescaling, we reduce the period  to one 
and make all coordinates and geometric parameters dimensionless. We also introduce
a domain $\Theta_\bullet \subset \bbR^3$, which is contained in the
cylinder $(-1/2, 1/2) \times (-\ell, \ell) \times \bbR$  
and has a smooth boundary $\partial 
\Theta_\bullet$ and compact closure $ \overline{\Theta_\bullet}
= \Theta_\bullet \cup \partial \Theta_\bullet$. We assume that the 
lower part $\Theta = \{ x \in \Theta_\bullet \, : \, z < 0 \}$ is non-empty
and denote its translates
\begin{eqnarray}
\Theta_j = \{ x \, : \, (y_1- j , y_2,z) \in \Theta \} \ \ \ \forall \, 
j \in \bbZ= \{0 , \pm1, \pm 2 , \ldots \} .  \label{1}
\end{eqnarray}
The periodic channel 
\begin{eqnarray}
\Pi^h = \Omega^h \cup {\textstyle \bigcup\limits_{j \in \bbZ}} \Theta_j 
\label{2}
\end{eqnarray}
consists, in addition to the deep pits \eqref{1}, of the thin periodic layer
\begin{eqnarray}
\Omega^h = \{ x = (y,z) \, : \,  0 > z > - h H(y), \ y \in \Omega^0\},
\label{3}
\end{eqnarray}
where $h>0$ is a small parameter. 
For the simplicity of the presentation (see Remark \ref{RR}) we assume that the subdomain
$\{x \in \Theta_\bullet \, : \, |z| < h_\bullet \}$ includes the straight
vertical cylinder  $\theta \times (- h_\bullet ,h_\bullet)$, where
$h_\bullet > 0$  is a constant and $\theta = \{ x \in \Theta_\bullet\, : \,
x_3 = z = 0\}$ is the cross-section of $\Theta_\bullet$, i.e., 
the surface $\Upsilon = \{x \in \partial \Theta
\, : \, |z| < h_\bullet \}$  is cylindrical and perpendicular to the plane $\{x\, :
\, z = 0\}$. Notice that we will not distinguish in the notation between the 
three-dimensional sets $\Omega^\circ, \theta, \partial \theta$ and their 
embeddings 
into the above-mentioned horizontal plane.

\begin{figure}
\begin{center}
\includegraphics[height=5cm,width=11cm]{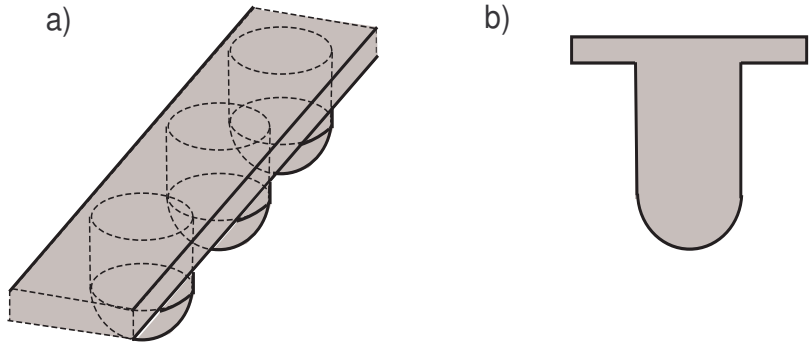}
\end{center}
\caption{a) Periodic channel, b) cross-section of a voluminous container, pothole}
\label{fig1}
\end{figure}

In the domain $\Pi^h$ we consider the water-wave problem, cf.\,\cite{KuMaVa}, which consists
of the Laplace equation for the velocity potential $u^h$,
\begin{eqnarray}
- \Delta u^h (x) = 0, \ \  x \in \Pi^h,   \label{u1}
\end{eqnarray} 
the Neumann (no penetration) boundary condition on the bottom and walls,
\begin{eqnarray}
\partial_\nu u^h(x) = 0, \ \ x \in \Sigma^h = \partial \Pi^h \setminus
\overline{\Omega^0}, \label{u2}
\end{eqnarray}
and the Steklov (kinematic) spectral condition on the free surface,
\begin{eqnarray}
\partial_z u^h(x) = \lambda^h u^h(x) , \ \ x =(y,0) \in \Omega^0. 
 \label{u3}
\end{eqnarray}
Here, $\Delta$ stands for the Laplacian, $\partial_\nu$ is the outward normal 
derivative, $\partial_z=\partial/\partial z$, and $\lambda^h = g^{-1}
k_h^2$ is a spectral parameter, where $g$ is the acceleration of 
gravity and $k_h$ is the physical wave number. 

As for the other general  notation used in this paper, we write 
$\bbN = \{1,2,3, \ldots \}$, $\bbN_0 = \{0 \} \cup \bbN$,
and  $\bbR_0^+$ for the set of non-negative real numbers. 
Given a domain $\Xi \subset \bbR^d$, the symbol $|\Xi|$ stands for its
volume in $\bbR^d$  
and   $(\cdot,\cdot)_\Xi$ stands for the natural scalar product in 
Lebesgue space $L^2(\Xi)$, and $H^k(\Xi)$, $k\in \mathbb{N}$,  for the 
standard Sobolev space of order $k$  on $\Xi$.  The norm of a function $f$ 
belonging to a Banach function space $X$ is denoted by $\Vert f ; X \Vert$. 
For $r >0$ and $a \in \mathbb{R}^N$, $B(a,r)$ (respectively, $S(a,r)$\,) 
stand  for the Euclidean ball (resp. ball surface) with centre $a$ and radius 
$r$. By $c, C$ (respectively, $c_k$, $C_k$, $c(k)$ etc.) we mean positive 
constants (resp. constants depending on a parameter $k$) which do not depend 
on  functions or variables appearing in the inequalities, but which may still 
vary from place to place.  The gradient and Laplace operators $\nabla$ and 
$\Delta$ (respectively, in $\nabla_y$, $\Delta_y $ etc.) act in the variable 
$x$ (respectively,  $y$ etc.).

\subsection{The goals of the paper} \label{sec1.2}
Our aim is to investigate the band-gap structure of the spectrum of the
water-wave problem \eqref{u1}--\eqref{u3} and in particular to detect
gaps in its essential spectrum. At the end of Section \ref{sec1} we will 
discuss the standard approach of the Floquet-Bloch-Gelfand-theory,
which transforms the problem into a model spectral problem in a bounded 
domain. The model problem will be treated by the methods of asymptotic
analysis in Section \ref{sec2}, where the asymptotic ans\"atze and their
essential terms are derived. In Section \ref{sec3}, we consider the 
relationship of the model and limit problems, where the latter 
corresponds to the case $h=0$. The main result, Theorem \ref{3AS} contains
the crucial estimate, up to corrections of order $O(h^2 (1 + |\ln h|))$, 
between the eigenvalues of the limit problem and those of the model problem 
(which  form the spectral bands of the original problem).  
Finally, in Section \ref{sec4} we consider the limit problem in more
detail under some additional assumptions. We
calculate its the eigenvalues and observe the existence of a gap, so that combining with Theorem \ref{3AS}, a gap
is opened also in the essential spectrum of the problem 
\eqref{u1}--\eqref{u3}.

Besides the Steklov spectral condition and periodic structure, one more
characteristic feature of the problem under consideration is 
described by junctions of massive bodies with thin ligaments. For example,
the dumbbell, which is a union of two massive domains connected by a thin
cylinder, is a classical object in asymptotic analysis, and the spectral
Laplace-Neumann problem and 
asymptotic expansions for eigenvalues and eigenfunctions have been 
considered in many papers.
Concerning asymptotic methods, we will here partly follow the approach in 
\cite{na526} and \cite{na576}, although the aims in these references are
different as they include topics like self-adjoint extension of the
problem operator and eigenvalue estimates also in the high-frequency
range, instead of spectral gaps corresponding to low frequencies in 
the present work. We mention the paper \cite{CCNT}, 
where the two-dimensional version of the present problem was treated
with result on the existence of arbitrarily many spectral gaps
for thin enough (small $h>0$) connecting canals of voluminous
containers. However, we emphasize that in the lower dimensional case the 
limit problem is an ordinary differential equation, the solutions of which
can be found  explicitly and quite precise information of its spectrum,
including spectral gaps, can be obtained by a quite elementary approach. 
Here, the limit problem is still a spectral Laplace problem, and precise
enough information on its spectrum is more difficult to gain. Also, it is
clear that the structure of the asymptotic terms is essentially different
and more complicated in the case $d=3$. 

Of course, neither the present paper nor \cite{CCNT} are   the first 
studies of spectral gaps by means of asymptotic analysis. Let us mention
\cite{na453,na530,BoPa,BkhMEP,NaOrPe}, where the detection of open gaps is based
on a periodic perforation of strips and cylinders or singular perturbations
of a similar type, an approach which will also be used in Section \ref{sec4}.

\subsection{Model problem in the periodicity cell} \label{sec1.3}
Our approach to the problem \eqref{u1}--\eqref{u3} is based on the 
Floquet-Bloch-Gelfand- (FBG-)transform, which converts the problem on the
periodic channel $\Pi^h$ into one on the periodicity cell   
\begin{eqnarray}
\varpi^h = \{ x \in \Pi^h \, : \, |y_1| < 1/2 \}. \label{01}
\end{eqnarray}
The boundary $\partial \varpi^h$ includes the free water 
surface $\omega^0 =  \{ x \in \Omega^0 \, : \, |y_1| < 1/2 \}$ and   
the lateral surface $\varsigma^h =  \{ x \in \Sigma^h \, : \, |y_1| < 1/2 \}$.
We recall that the FBG-transform is defined by 
\begin{eqnarray*} 
u(x)\mapsto\ U(x,\eta)=\frac{1}{\sqrt{2\pi}}\sum_{j\in \mathbb{Z}}
e^{-i\eta j} u(x_1+j, x_2,x_3),
\end{eqnarray*} 
where $x = (x_1, x_2,x_3)\in \Pi^{h}$ on the left while,  on the right, we have 
$x \in\varpi^h$ and $\eta\in[-\pi, \pi]$; this is  the dual variable or 
the Floquet parameter. For more details, see \cite{Gel} and e.g. 
\cite[XII.16]{RS}, \cite[\S\, 3.4]{NaPl}, \cite[Cor. 3.4.3]{NaSpec},  
\cite[Sec.\,2.2]{Kuchbook}. 
Applying the FBG-transform to the problem (\ref{u1})--(\ref{u3}), we obtain a 
family of model problems in the periodicity cell $\varpi^h$ parametrized by 
$\eta$, 
\begin{eqnarray}
-\Delta U^h(\eta; x) &=& 0, \quad x\in \varpi^h, 
\label{7}  \\
\partial_{\nu}  U^h(\eta; x) &=& 0, \quad x\in \varsigma^{h},
\label{8} \\
\partial_z  U^h(\eta; x)&=& \Lambda^h(\eta)  U^h(\eta; x), 
\quad x\in \omega^0 , 
\label{9} \\
 U^h(\eta; x) \Big\vert_{y_1 =1/2}& =& e^{i\eta}   U^h(\eta; x)
 \Big\vert_{y_1 = - 1/2} ,            
\label{10} \\
\frac{\partial U^h}{\partial x_1}  (\eta; x) \Big\vert_{y_1 =1/2}
&=& e^{i \eta}
\frac{\partial U^h}{\partial x_1}  (\eta; x) \Big\vert_{y_1 = - 1/2} .  \label{11}
\end{eqnarray}
Here, 
$\Lambda =\Lambda(h,\eta) $ is a new notation 
for the spectral parameter $\lambda$.

As is well known, the FBG-transform 
establishes an isometric isomorphism 
\begin{equation*}
L^2(\Pi^{h})\simeq L^2(0,2\pi; L^2(\varpi^h)),
\end{equation*}
where $L^2(0,2\pi;B)$ is the Lebesgue space of functions with values in the Banach space $B$,  endowed with the norm
\begin{equation*}
\| f ;L^2(0,2\pi; B)\|=\Big(\int_0^{2\pi}\|f (\eta); B \|^2d\eta
\Big)^{1/2}\,.
\end{equation*}
We denote 
by $H_\eta^1(\varpi^h)$, where $\eta \in [0, 2\pi)$, the subspace of the Sobolev space  
$H^1(\varpi^h)$
consisting of functions satisfying the quasiperiodic boundary condition 
\eqref{10}. 
The FBG-transform is also an isomorphism from the Sobolev space 
$H^1(\Pi^{h})$ onto \hfill\break 
$H^1(0,2\pi; H^1_{\eta}(\varpi^h))$;
see the references above.  

Our approach to the spectral properties of the model  is 
similar to  \cite{na526,NaTa} and others. Given $\eta \in [\pi, \pi]$, 
we write the variational formulation of the problem \eqref{7}--\eqref{11}
for the unknown function $U \in H_\eta^1(\Omega^h)$ and the spectral
parameter $\Lambda(\eta)$  as
\begin{eqnarray}
\label{15}
& & (\nabla U^h(\eta;  \cdot)   , \nabla \psi^h(\eta ;  \cdot) )_{\varpi^h} 
\nonumber \\
&=& \Lambda^h(\eta) 
(U^h(\eta ;  \cdot) ,\psi^h(\eta ;  \cdot) )_{\omega^0} \, , \  
\forall \, \psi^h(\eta;  \cdot) \in H^1_{\eta}(\varpi^h) .
\end{eqnarray}
We denote in the sequel by $\cH_\eta^h$ the space $H_{\eta}^1 (\varpi^h)$ endowed with  the new scalar product 
\begin{equation}
\langle u^h, \psi^h \rangle_{h} = 
\big (\nabla u^h, \nabla \psi^h\big)_{\varpi^h} 
+ h \big(u^h, \psi^h \big)_{\omega^0} , \label{58}
\end{equation} 
where $\eta$ is omitted from the notation and the inner product in 
$\omega^0$ is understood in the sense of traces, 
and define  a self-adjoint, positive  operator 
$\cK_\eta^h : \cH_\eta^h \to \cH_\eta^h $ by using
the identity 
\begin{equation}
\langle \cK_\eta^h u^h, \psi^h\rangle_{h} = 
\big(u^h,\psi^h\big)_{\omega^0} \ \ \forall \, 
u^h, \psi^h \in  \cH_\eta^h. \label{59}
\end{equation}
The operator $\cK_\eta^h$ is also compact due to the compactness
of the embedding $H^1(\varpi^h) \subset L^2(\omega^0)$. 
The integral identity \eqref{15} corresponding to the problem 
\eqref{7}--\eqref{11} is then equivalent with the standard spectral 
equation in the Hilbert space $\cH_\eta^h$,
\begin{eqnarray*} 
\cK_\eta^h U^h = \kappa^h U^h , 
\end{eqnarray*} 
with the new spectral parameter 
\begin{equation}
\kappa^h = \kappa^h (\eta)  
= \big( \Lambda^h(\eta) + h \big)^{-1}.\label{60}
\end{equation}
According to \cite[Thm.\,10.1.5,\,10.2.2]{BiSo}  the spectrum of 
$\cK_\eta^h$ consists of null, which is the only point 
in the essential spectrum,  and a positive sequence of eigenvalues 
belonging to  the discrete spectrum
\begin{eqnarray}
\kappa_1^h \geq \kappa_2^h \geq \ldots \geq  \kappa_m^h
\geq \ldots \to + 0 .\label{61}
\end{eqnarray}
Formula \eqref{60} and the properties of the sequence 
$\big(\kappa_k^h (\eta) \big)_{k=1}^\infty$ mean that the eigenvalues 
of the problem \eqref{7}--\eqref{11} form an unbounded sequence
\begin{equation}
\label{12}
0 \leq \Lambda_1^{h} (\eta) \le \Lambda_2^{h}(\eta)  \le \ldots\le\Lambda_k^{h} (\eta) 
\le\ldots\to+\infty ,
\end{equation}
where multiplicities have been taken into account. We denote for every $k$  by 
$ U_k^h(\eta; \cdot)   \in H_\eta^1(\varpi^h)$ the 
eigenfunction  corresponding to $\Lambda_k^h (\eta) $, which fulfil the
normalization and orthogonality conditions
\begin{eqnarray}
\big(   U_m^h(\eta; \cdot) ,  U_n^h(\eta; \cdot) \big)_{\omega^0}
= \delta_{m,n}  , \ \ m,n\in\bbN,\label{16}
\end{eqnarray}
where on the right-hand side there is the Kronecker delta.

Finally, the functions $\eta\mapsto\Lambda_k^{h}(\eta)$  are continuous and 
$2\pi$-periodic (see for example \cite[Ch.\,9]{Kato}, 
\cite[Sec.\,3.1]{Kuchbook}). Hence, the spectral
bands 
$
{\sf B}_m^{h}=\{\Lambda_m^{h}(\eta):\eta\in [0,2\pi)\}
$
of the problem \eqref{u1}--\eqref{u3} indeed are compact intervals,
and, by well-known principles of the FBG-theory, the essential spectrum
$S_{\rm ess}^h$ of the problem \eqref{u1}--\eqref{u3} can  be presented as
\begin{eqnarray*}
S_{\rm ess}^h = {\textstyle \bigcup\limits_{m \in \bbN} }{\sf B}_m^{h}.
\end{eqnarray*}

\begin{figure}
\begin{center}
\includegraphics[ height=5cm,width=12cm]{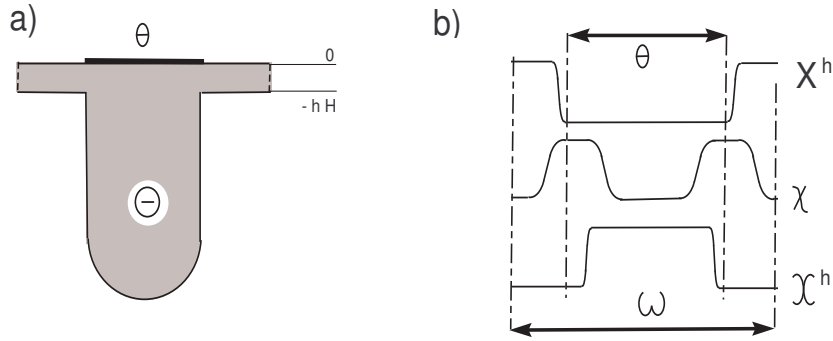}
\end{center}
\caption{a) Cross-section of the structure of a pothole, b) graphs of cut-off functions.}
\label{fig2}
\end{figure}

\section{Formal asymptotic analysis}  \label{sec2}
\subsection{Asymptotic ans\"atze} \label{sec2.1}

We search for an eigenvalue of the problem \eqref{7}--\eqref{11} in the form
\begin{eqnarray}
\Lambda^h(\eta) = h \mu (\eta) +\ldots \ .\label{21}
\end{eqnarray}
Here, and later in this section, dots stand for higher order terms, which
are inessential for our preliminary analysis. Inside the thin part 
$\varpi_\sharp^h = \varpi^h \setminus \overline \Theta$ of the periodicity
cell, we employ the standard asymptotic decomposition in shallow
water domains
\begin{eqnarray}
U^h(\eta;x) = v^0(\eta;y) + h^2 v'(\eta ; y, \zeta) + \ldots
\ \ \ \mbox{in} \ \varpi_\sharp^h, \label{22}
\end{eqnarray}
cf.\,\cite{na470,na526,na576,CCNT} and others. Here, $\zeta = h^{-1} z$ 
is the stretched vertical coordinate. We note that
\begin{eqnarray}
\partial_\nu = \big( 1 + h^2 |\nabla_y H(y)|^2 \big)^{-1/2}
\big( -\partial_z - h \nabla_y H(y) \cdot \nabla_y \big)  \label{23}
\end{eqnarray}
at the bottom part $\varsigma_\sharp^h = \{ x \in \omega^0 \setminus
\overline{\theta} \, : \, z = - hH(y) \}$ and insert the ans\"atze into the 
equation \eqref{7}, restricted to $\varpi_\sharp^h$, and to the 
boundary conditions \eqref{8} and \eqref{9}, restricted to
$\varsigma_\sharp^h$ and $\omega_\sharp^0 = \omega^0 \setminus 
\overline \theta $, respectively. Separating terms of the same order in
$h$ yields the problems  
\begin{eqnarray*}
- \partial_\zeta^2 v^0(\eta;y) = 0 \ \mbox{for} \  \zeta \in (-H(y),0), \ \
\  \partial_\zeta v^0(\eta;y) = 0 \ \mbox{for} \  \zeta = -H(y),0 
\end{eqnarray*}
(which trivially holds, since $v$  is independent of $\zeta$ in \eqref{22}),
and
\begin{eqnarray}
& & -  \partial_\zeta^2 v'(\eta;y,\zeta) = \Delta_y v^0(\eta ;y)  ,   \  \zeta \in (-H(y),0), \nonumber \\
& &  - \partial_\zeta  v'(\eta;y, - H(y) ) =  \nabla_y H(y) \cdot \nabla_y v^0(\eta;y) , \
- \partial_\zeta  v'(\eta;y, 0 ) = \mu v^0(\eta ;y) .
\label{24}
\end{eqnarray}
Equations \eqref{24} form a Neumann problem for an ordinary differential
equation, and it has a solution, if and only if the compatibility condition
\begin{eqnarray*}
\int\limits_{-H(y)}^0 \Delta_y v^0(\eta;y) d\zeta + \nabla_y H(y) \cdot\nabla_y v^0(\eta;y) + \mu(\eta) v^0(\eta;y) = 0 
\end{eqnarray*}
is satisfied. This can be written as the differential equation
\begin{eqnarray}
- \nabla_y \cdot H(y)\nabla_y v^0(\eta;y) = \mu(\eta) v^0(\eta;y),
\ \ \ \ y \in \omega_\sharp^0. 
\label{25}
\end{eqnarray}
Due to the form of the ansatz \eqref{22}, the conditions 
\eqref{8}, \eqref{10} and \eqref{11} lead to the following boundary
and quasiperiodicity conditions
\begin{eqnarray}
& & \pm \frac{\partial v^0}{\partial y_2}(\eta; y_1 , \pm \ell) = 0, \ \ \ |y_1| < 1/2 , \label{260} 
\\
& & v^0 \Big(\eta; \frac12 ,y_2 \Big) = e^{i \eta} v^0 
\Big(\eta;- \frac12 ,y_2 \Big) , \ \ \ 
|y_2| < \ell \label{261}
\\
& & \frac{\partial v^0}{\partial y_1} \Big(\eta; \frac12,y_2 \Big) = e^{i \eta} 
\frac{\partial v^0}{\partial y_1}\Big(\eta;-\frac12 ,y_2 \Big), \ \ \ 
|y_2| < \ell . \label{26}
\end{eqnarray}

It remains to derive a  condition on the interior boundary $\partial \theta$
of the domain $\omega_\sharp^h \subset \bbR^2$. To this end we use the 
asymptotic ansatz, cf. \cite{na526,na576,CCNT}, 
\begin{eqnarray}
U^h(\eta;x)= V^0(\eta;x) + h V'(\eta;x) + \ldots \label{27}
\end{eqnarray}
in the massive part  $\Theta$ of the periodicity cell \eqref{01}.
We insert the ans\"atze \eqref{27} and \eqref{21} into the equation
\eqref{7} in $\Theta$, into the Steklov condition \eqref{9} in $\theta$
and into the Neumann condition \eqref{8} on the surface $\varsigma^h \cap\partial \Theta$. The curved two-dimensional ring
\begin{eqnarray}
\upsilon^h = \{ x \, : \, y = \partial \theta, -h H(y)< z < 0 \} 
\label{03}
\end{eqnarray}
does not appear in the boundary condition \eqref{8} and it also disappears
on the limit $h \to +0$ so that it is natural to compose the problems
\begin{eqnarray}
- \Delta V^0(\eta; x ) = 0 \ \mbox{for} \ x \in \Theta, \ \
\  \partial_\nu V^0(\eta;x) = 0 \ \mbox{for} \  x \in \partial \Theta
\setminus \partial \theta ,
\label{28}
\end{eqnarray}
and
\begin{eqnarray}
& & - \Delta V'(\eta; x ) = 0 \ \mbox{for} \ x \in \Theta, \ \
\  \partial_\nu V'(\eta;x) = 0 \ \mbox{for} \  x \in \partial \Theta
\setminus \overline \theta
\nonumber \\
& & - \partial_z  V'(\eta; x ) = \mu(\eta) V^0(\eta; x) , \ \ x \in \theta .
\label{29}
\end{eqnarray}

According to \eqref{28}, the main asymptotic term in \eqref{27} 
is a constant,
\begin{eqnarray}
V^0(\eta;x) = a^0(\eta); \label{30}
\end{eqnarray}
however, the Neumann problem \eqref{29} has no bounded solution, if
$\mu(\eta) \not= 0$ and $a^0(\eta) \not=0$. This leads to the conclusion that 
the problem \eqref{29} has to be modified (remember how we treated 
the ring \eqref{03}).

\subsection{Boundary layer.}  \label{sec2.2}
\begin{figure}
\begin{center}
\includegraphics[width=10cm]{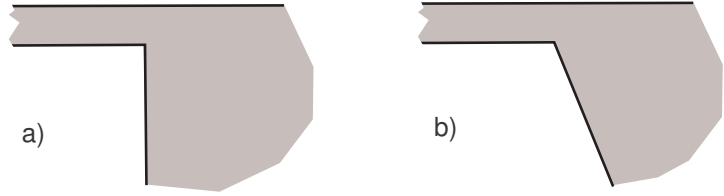}
\end{center}
\caption{a) The scaled domain $\Xi$, b) non-perpendicular case.}
\label{fig3}
\end{figure}

We will link the outer expansion \eqref{22} and \eqref{27} by using the
method of matched asymptotic expansions, see the monographs \cite{VM,Ilin} and
the papers \cite{na526,na576,CCNT}, where applications to the water-wave
problem are considered. To construct the inner expansion in the vicinity of
$v^h$, we introduce the local curvilinear coordinates $(s,n,z)$ in a 
three-dimensional neighborhood $\cV$ of the contour $\partial \theta$, such that
$s$ is the arc length on $\partial \theta$, $n$ is the oriented distance to
$\partial \theta$ on $\omega^0$ with $n > 0$ inside $\theta$, and $z = x_3$.
We keep $s$ unscaled but change $n$ and $z$ into the stretched coordinates
\begin{eqnarray}
\xi = (\xi_1, \xi_2) = \frac{1}{h} \frac{1}{H(s)} (n,z) \ \ \ 
\mbox{with}  \ H(s) = H(y)\big|_{\partial\theta } . \label{31}
\end{eqnarray}
The coordinate change $x \mapsto (s, \xi) $ and the formal replacement 
$h= 0$ transform the subdomain $\varpi^h \cap \cV$ into the set
\begin{eqnarray*}
\big\{ (s,\xi) \, : \, s \in \partial \theta, \xi \in \Xi \big\}
\end{eqnarray*}
where $\Xi$ is the union of the fourth quadrant of the plane,
$\bbK = \{ \xi \in \bbR^2 \, : \, \xi_1 > 0, \xi_2 <0 \}$ and the
semi-infinite strip $\bbP= \{ \xi \in \bbR^2 \, : \, \xi_1 \leq 0, \xi_2 
\in (-1,0)  \}$ of unit width (recall the denominator $H(s)$ in \eqref{31};
see Fig.\,\ref{fig3},a)).

We have $\partial_\nu = \partial_n = h^{-1} \partial / \partial  \xi_1$ for
$y \in \partial \theta$, $z \in (-h_\bullet,0)$ and $\partial_z = h^{-1} 
\partial / \partial \xi_2$ on $\omega^0$. The Laplacian $\Delta_x$ is written
in the curvilinear coordinates as
\begin{eqnarray}
\big( 1 + n \varkappa (s) \big)^{-1} \Big( \frac{\partial}{\partial s}
\big( 1 + n \varkappa (s) \big)\frac{\partial}{\partial s} + 
\frac{\partial}{\partial n}    \big( 1 + n \varkappa (s) \big)^{-1}
\frac{\partial}{\partial n} \Big) + \frac{\partial^2}{\partial z^2}  , \label{313}
\end{eqnarray}
where $\varkappa(s)$ is the curvature of the contour $\partial \theta \ni s$. 
Thus, the coordinate change $x \mapsto \xi$ together with formulas \eqref{313} and 
\eqref{23} and passing formally to the limit $h \to +0$ turn \eqref{7}--\eqref{9} 
into the Neumann problem\begin{eqnarray}
- \Delta_\xi w(\xi) = 0, \ \ \xi \in \Xi, \ \ \ 
\partial_{\nu(\xi)} w(\xi) =0, \ \ \xi \in \partial \Xi \setminus P, \label{32}
\end{eqnarray}
where $P= (0,-1)$; this is the only corner point with opening $\alpha = 3\pi/2$
on the boundary $\partial \Xi$. 

\begin{remark}
\label{RR} If the surface $\partial \Theta$ is not perpendicular to the plane
$\{x \, : \, z = 0\}$ on the contour $\partial \theta$, then, instead of the 
union $\Xi_{ 3 \pi /2} = \bbK \cup \bbP$ (Fig.\,\ref{fig3},b)) we obtain the domain 
\begin{eqnarray*}
\Xi_{\alpha(s)} = \bbK_{\alpha(s)} \cup \{ \xi \, : \,
\xi_1 \in \bbR, \xi_2 \in (-1,0) \}, 
\end{eqnarray*}
where $\bbK_{\alpha(s)} = \{ \xi \, : \, \xi_2 < 0, \xi_1 > - \xi_2
\cot (\alpha(s)) \}$ is the angle of opening $\alpha(s) \in (\pi, 2\pi)$. 
In this case the formal asymptotic analysis would remain almost the same, but the
justification procedure of Section \ref{sec3} would become more involved 
with more cumbersome calculations. 
\end{remark}

To follow the matching procedure in \cite{na526,CCNT}, we need to find
the solutions of the problem \eqref{32} with linear behavior at the infinity
in $\bbP$. One of them is obvious, the constant function. However,  there also exists a harmonic function $W$ in  $\Xi$ satisfying the Neumann condition such that
\begin{equation}
W(\xi) = \left\{ 
\begin{array}{ll}
\xi_1 + c_\Xi + O(e^{\pi \xi_1} ) , & \xi_1 \to - \infty, \ \xi \in \bbP \\
2 \pi^{ -1}  \ln |\xi| +  O(|\xi|^{-1}) , & |\xi_1| \to + \infty, \ \xi \in \bbK ,
\end{array} \right.  
\label{33}
\end{equation}
see for example \cite[\S\,2]{na576}.
Here $c_\Xi$ is an absolute constant which can be found by solving the problem
by an appropriate conformal mapping, but we do not need the exact value in the
following. Note that the solution is made unique by the requirement that the 
constant term of the solution, which remains undetermined by the Neumann 
condition,  equals 0 in the quadrant $\bbK$. 

Furthermore, it follows from the Kondratiev theory (see \cite{Ko} and, e.g., \cite[Ch.\,2]{NaPl}) that all solutions of
\eqref{32} with at most polynomial growth rate in $\bbP$ and logarithmic in 
$\bbK$ are given by the linear combinations
\begin{eqnarray*}
c_0 + c_1 W(\xi), \ \ \ c_0,c_1 \in \bbR.
\end{eqnarray*}

Let us now realize the matching procedure. The Taylor formula in the variable $n$
converts the outer expansion \eqref{22} in the thin domain $\varpi_\sharp^h$ into
\begin{eqnarray}
U^h(\eta; x) &= & v(\eta; s,0) + n \partial_n v(\eta;s,0) + \ldots +
h^2v'(\eta; s,0) + \ldots 
\roweq
v(\eta;s,0) + hH(s) \partial_n v(\eta; s,0)+ \ldots ,\label{34} 
\end{eqnarray}
where $v(\eta; s,n)$ is the function $v(\eta; \cdot)$ written in the curvilinear
coordinates. Comparing terms of order $1=h^0$ in the outer expansion \eqref{34} 
and \eqref{27} with the inner expansions
\begin{eqnarray}
U^h (x) = c^0(\eta ; s) + h c^1(\eta;s) 
W\big(h^{-1}H(s)^{-1} n, (h^{-1}H(s)^{-1}z \big) + \ldots , \label{35}
\end{eqnarray}
we conclude that 
\begin{eqnarray}
v^0(\eta;s,0) = c^0(\eta;s) = a^0 (\eta), \label{36}
\end{eqnarray}
where $a^0(\eta)$ is taken from \eqref{30}.

To continue, we compare \eqref{34} and \eqref{35} at the level $h$. 
The asymptotics \eqref{33} in $\bbP$ requires that 
\begin{eqnarray*}
c^1(\eta; s) = H(s) \partial_nv^0(\eta; s,0).
\end{eqnarray*}
Hence, in view of the representation \eqref{33} in $\bbK$, the correction term
of the outer expansion \eqref{27} in the massive part $\Theta$ of
the periodicity cell $\varpi^h$ must have the following behavior at edge
$\partial \theta$ of the boundary $\partial \Theta$:
\begin{eqnarray}
V'(\eta;x) = \frac{2}{\pi} H(s) \partial_n v^0(\eta;s,0) \ln r +
O(1)  \ \ \ \mbox{as} \ r = (n^2 +z^2)^{1/2} 
\to + 0 .    \label{37}
\end{eqnarray}

It should be mentioned that the thin ring \eqref{03} does not involve any
boundary condition and it shrinks into the edge $\partial \theta$ of the 
container $\Theta$ so that it is not possible to make a priori any conclusion
on the behavior of $V'( \eta ; x)$ as $x$ approaches $\partial \theta$. 
At the same time, the matching procedure leads to the relation \eqref{37},
which completes the formulation of the Neumann problem for $V'$. 

A solution of the problem \eqref{29}, \eqref{37} can be found in the form
\begin{eqnarray*}
V'(\eta;x) = \chi(x) \frac{2}{\pi} H(s) \partial_n v^0(\eta;s,0) \ln r +
\widehat V'(\eta;x) ,
\end{eqnarray*}
where $\chi$ is a smooth cut-off function which is equal to 1 in the
vicinity of $\partial \theta$ and has support in $\overline \Theta \cap \cV$.
The function  $\widehat V'(\eta,x)$ can be found from the Poisson equation with
the right-hand side $\widehat f$ such that $r^\delta \widehat f \in L^2(\Theta)$
and with Neumann data in $L^2(\partial \theta)$. Here, $r(x)  ={\rm dist}\,
(x, \partial \theta)$ and $\delta > 0$  is arbitrary. Hence, the existence of 
$\widehat V' \in H^1 (\Theta)$ follows by posing one compatibility condition. 
This can of course be derived from the original problem \eqref{29}, \eqref{37} by 
integrating by parts in the truncated domain 
$\Theta (\varrho) = \{ x \in \Theta \, : \, r(x) > \varrho \}$ and sending
$\varrho $ to $+0$. Indeed,
\begin{eqnarray*}
\mu \, {\rm mes}_2 (\theta) \,  a^0(\eta) &=& 
\lim\limits_{\varrho \to +0} \int\limits_{ \{y\in \theta : r > \varrho\} } 
\partial_z V'(\eta; y,0) \, dy
\roweq
\lim\limits_{\varrho \to +0} \int\limits_{\{ x \in \Theta \cap \cV: r = 
\varrho\} } 
\partial_r V'(\eta;x) \, ds_x
\roweq
\int\limits_{\partial \theta} H(s) \partial_n v^0(\eta;s,0) ds \frac{2}{\pi}
\int\limits_0^{\pi/2} d\varphi = 
\int\limits_{\partial\theta}  H(s) \partial_n v^0(\eta;s,0) ds ,
\end{eqnarray*}
where $\varphi$ is the angular variable of the polar coordinate system 
$( r, \varphi)$ in planes perpendicular to the contour $\partial \theta$. 

As a result of our calculations, we conclude by associating 
the following conditions on $\partial \theta$ to the problems 
\eqref{25}--\eqref{26} ,
\begin{eqnarray}
& & v^0(\eta;s,0) = a^0(\eta) , \ \ \ s \in \partial \theta, 
\row
\int\limits_{\partial \theta} H(s) \partial_n v^0(\eta; s,0) ds 
=  \mu(\eta) a^0(\eta) {\rm mes}_2(\theta) . 
\label{38}
\end{eqnarray} 
The quantity $a^0(\eta)$ has not been fixed yet, but it must be found
by solving the obtained eigenvalue problem \eqref{25}--\eqref{26}, \eqref{38},
which we call the limit problem.

\subsection{The spectrum of the limit problem.} \label{sec2.3}
To present the variational formulation of the problem 
\eqref{25}--\eqref{26}, \eqref{38} we introduce the following subspace of the
Sobolev space $H^1 (\omega_\sharp^0)$:
\begin{eqnarray}
\ \ \ 
\cH(\eta) = \big\{ v^0 \in H^1(\omega_\sharp^0) &:&
v^0 \mbox{ satisfies \eqref{261} and
 is constant on } \partial \theta \subset 
\partial \omega_\sharp^0 \big\}. \label{39}
\end{eqnarray}
We denote the a priori unknown constant in \eqref{39} by $v_\theta^0 \in \bbC$.
Owing to \eqref{38} and \eqref{39} we have
\begin{eqnarray*}
\int\limits_{\partial \theta} H(s) v^0 (\eta; s,0) \overline{ \psi
(\eta; s,0)} ds = \mu(\eta) v_\theta^0 \overline{\psi_\theta} 
{\rm mes}_2(\theta) 
\end{eqnarray*}
for any solution $v^0 \in H^2(\omega_\sharp^0)$ of the limit problem
and test function $\psi (\eta; \cdot) \in \cH(\eta)$, hence, multiplying
\eqref{25}  by $\overline{ \psi(\eta;x) }$, integrating by parts, and taking
into account the boundary and quasiperiodicity conditions \eqref{25}--\eqref{26}
yield the integral identity
\begin{eqnarray}
& & \big( H \nabla_y v^0(\eta; \cdot) , \nabla_y \psi(\eta; \cdot) 
\big)_{\omega_\sharp^0} 
\roweq
\mu(\eta) \big( (v^0(\eta ; \cdot) ,  \psi(\eta;\cdot) \big)_{\omega_\sharp^0}
+ v_\theta^0(\eta) \overline{ \psi_\theta(\eta)} {\rm mes}_2(\theta)
\ \ \ \forall \, \psi (\eta; \cdot) \in \cH(\eta) ,  \label{40}
\end{eqnarray}
where $( \cdot, \cdot)_{\omega_\sharp^0}$ is the natural scalar product of the
Lebesgue space $L^2({\omega_\sharp^0})$.

Since the embeddings of $H^1( {\omega_\sharp^0})$ into $L^2({\omega_\sharp^0})$ and
$L^2(\partial \theta)$ are compact, we observe that \eqref{39} is a closed subspace
of $H^1({\omega_\sharp^0})$. Moreover, the sesquilinear form on left-hand side of 
\eqref{40}  is positive, thus Hermitian,  and closed in $H^1({\omega_\sharp^0})$.
These observations obviously imply the following assertion.

\begin{proposition} \label{P1}
The spectrum of the problem \eqref{40}  is discrete and it consists of the
positive monotone unbounded sequence of eigenvalues
\begin{eqnarray}
0 \leq \mu_1(\eta) \leq \mu_2(\eta)\leq \ldots \leq \mu_m(\eta)  
\leq \ldots \to +\infty  \label{41}
\end{eqnarray}
where multiplicities are taken into account. The corresponding eigenfunctions
$v_1^0(\eta; \cdot)$, $v_2^0(\eta; \cdot)$, $\ldots , v_m^0(\eta; \cdot), 
\ldots \in \cH(\eta)$ can be subject to normalization and orthogonality conditions
\begin{eqnarray}
\big( v_m^0(\eta; \cdot) , v_n^0(\eta; \cdot) \big)_{\omega_\sharp^0}
+ v_{m\theta}^0(\eta) \overline{ v_{n \theta}^0(\eta) }{\rm mes}_2(\theta) 
= \delta_{m,n} \ \ \ \ m,n \in \bbN .  \label{43}
\end{eqnarray}
\end{proposition}

Since the boundary $\partial \theta$ was assumed smooth, also the eigenfunctions
$v_m^0(\eta, \cdot)$ are smooth, see \cite{ADN}, and therefore the eigenpairs
$\{ \mu_m(\eta), v_m^0(\eta ,\cdot)\} $ of the variational problem \eqref{40}
are solutions of the differential problem \eqref{25}--\eqref{26}, \eqref{28}, too. 

The first eigenvalue and eigenfunction of the limit problem satisfy the relations
\begin{eqnarray*}
0 = \mu_1(0) < \mu_1(\eta) \ \mbox{for} \ \eta\in [-\pi,0) \cup
(0, \pi], \ \ \ \ v_1^0(0;y) = (2 \ell)^{-1/2} . 
\end{eqnarray*}
It is remarkable that the first eigenpair $\{ \Lambda_1^h(0) , U_1^h(0,x)\}$of the original problem \eqref{7}--\eqref{11} with $\eta=0$ has the same form
$\{0, (2 \ell)^{-1/2}\}$.

\section{Justification of asymptotics} \label{sec3}
\subsection{Convergence theorem} \label{sec3.1}

We consider an eigenpair $\{ \Lambda_n^h(\eta) , U_n^h(\eta; \cdot ) \}$ of 
the problem \eqref{15}, where  $U_n^h(\eta ; \cdot)$ is normalized as in
\eqref{16} and 
\begin{eqnarray}
\Lambda_n^h(\eta) \leq C^{(n)} h \ \ \ \mbox{for } h \in (0, h^{(n)}]
\label{49}
\end{eqnarray}
for some positive $h^{(n)}$ and $C^{(n)}$; the bound \eqref{49} will be 
derived in Remark \ref{remNEW}. Hence, there exists a positive sequence $\{h_j\}_{j \in \bbN}$ tending to 0
such that
\begin{eqnarray}
h_j^{-1} \Lambda_n^{h_j} (\eta) \to \widehat{\mu}_n (\eta) \in [0, + \infty).
\label{50}
\end{eqnarray}
In what follows we omit the argument $h$ and the index $n$ from the notation and also write
simply $h$ instead of $h_j$. 

We pick up a function $\psi \in C^\infty (\overline{\omega_\sharp^0})$ which satisfies
the first quasiperiodicity condition \eqref{261} and coincides with $\psi_\theta \in \bbC$ on 
$\partial \theta$. We insert the special test function 
\begin{equation*}
\psi^h(x) = \left\{ 
\begin{array}{ll}
\psi(y) , & \ \ x \in \varpi_\sharp^h \\
\psi_\theta, & \ \ x \in \Theta ,
\end{array}
\right. 
\end{equation*}
into the integral identity \eqref{15} and rewrite it as follows:
\begin{eqnarray}
& & h^{-1} \int\limits_{\varpi_\sharp^h} \nabla_y U^h(x) \cdot \overline{ \nabla_y 
\psi(y) } dx - h^{-1} \Lambda^h
\int\limits_{\omega_\sharp^0}  U^h(y,0) \overline{\psi(y) } dy 
\roweq
h^{-1} \Lambda^h \int\limits_{\theta}  U^h(y,0)  \overline{\psi(y) } dy  .\label{52}
\end{eqnarray}
We set $\Psi^h = U^h$ in  \eqref{15} and use \eqref{49}, \eqref{16} to  derive the estimate
\begin{eqnarray*}
\Vert \nabla U^h ; L^2(\varpi^h) \Vert^2 = \Lambda^h \Vert U^h ; L^2(\omega^0) \Vert^2
= \Lambda^h \leq c h 
\end{eqnarray*}
so that the Poincar\'e and trace inequalities in the massive part $\Theta$ yield
\begin{eqnarray}
& & U^h(x) = U_0^h + U_\perp^h(x) , \ \ U_0^h \in \bbC, \ \ 
\int\limits_\Theta U_\perp^h(x) dx = 0, 
\row
\Vert U_\perp^h; L^2(\Theta) \Vert^2 \leq c \Vert \nabla U_\perp^h; L^2(\Theta) \Vert^2
= c \Vert \nabla U^h; L^2(\Theta) \Vert^2 
\rowleq
c \Lambda^h \Vert U^h; L^2(\omega^0) \Vert^2 \leq ch, 
\row
\Vert U_\perp^h; L^2(\theta) \Vert^2 \leq c \Vert U_\perp^h; H^1(\Theta) \Vert^2
\leq ch . \label{521}
\end{eqnarray}
Hence, we can pass to a subsequence of $\{ h_j\}_ {j \in \bbN}$, and still 
keep the above convention on the notation, so that
\begin{eqnarray}
U^h\big|_ {z = 0} \to \widehat v_\theta \in \bbC \ \ \ \mbox{strongly in } L^2(\theta).
\label{53}
\end{eqnarray}
Therefore we obtain for the right-hand side $I_{\rm ri}^h(\psi)$ of \eqref{52}
\begin{eqnarray}
I_{\rm ri}^h(\psi) \to \widehat \mu {\rm mes}_2 (\theta) \widehat v_\theta 
\overline{\psi_\theta}. \label{54}
\end{eqnarray}

To process the left-hand side $I_{\rm le}^h(\psi)$ of \eqref{52} we define the 
stretched domain $\varpi_\sharp^1 = \{ (y, \zeta) \, : \, 
y \in \omega_\sharp^0, - H(y) < \zeta < 0 \}$ and observe that the function
$\bfU^h(y ,\zeta) = U^h(y, h^{-1} z)$ satisfies
\begin{eqnarray}
& & \Vert \bfU^h(\cdot, 0) ; L^2(\omega_\sharp^0) \Vert^2 
+ \Vert \nabla_y \bfU^h(\cdot, 0) ; L^2(\varpi_\sharp^1) \Vert^2 +
h^{-2} \Vert \partial_\zeta \bfU^h(\cdot, 0) ; L^2(\varpi_\sharp^1) \Vert^2 
\roweq
\Vert U^h(\cdot, 0) ; L^2(\omega_\sharp^0) \Vert^2 +
h^{-1} \Vert \nabla_y  U^h ; L^2(\varpi_\sharp^h) \Vert^2 
+ h^{-1} \Vert \partial_z  U^h ; L^2(\varpi_\sharp^h) \Vert^2 
\rowleq(1+ h^{-1} \Lambda^h ) \Vert U^h(\cdot, 0) ; L^2(\omega^0) \Vert^2 
\leq C. \label{04}
\end{eqnarray}
Thus, we again find a  subsequence  of $\{ h_j\}_ {j \in \bbN}$  such that
\begin{eqnarray}
& & \bfU^h \to \widehat v \ \ \ \mbox{weakly in } H^1(\varpi_\sharp^1) 
\row
\bfU^h(\cdot, 0)  \to \widehat v \ \ \ \mbox{strongly in } 
L^2(\omega_\sharp^0) ,
\label{55}
\end{eqnarray}
where  $\widehat v  \in  H^1(\varpi_\sharp^0)$ is a function independent of $\zeta$
(recall the factor $h^{-1}$ of the norm $\Vert \partial_\zeta \bfU^h ; L^2(
\varpi_\sharp^1)\Vert$ in \eqref{04}). Finally,
\begin{eqnarray}
I_{\rm le}^h(\psi) \to 
\int\limits_{\omega_\sharp^0} H(y) \nabla_y \widehat v (y) \cdot \nabla_y \psi(y) dy - 
\widehat \mu \int\limits_{\omega_\sharp^0} H(y)  \widehat v (y) \cdot \psi(y) dy 
\label{56}
\end{eqnarray}
and, moreover, from the normalization of $U^h$ and the strong convergence in 
\eqref{53} and \eqref{55} we derive the equality
\begin{eqnarray}
\Vert \widehat v ; L^2(\omega_\sharp^0) \Vert^2 + |\widehat v_\theta |^2
{\rm mes}_2 (\theta) = 1. \label{57}
\end{eqnarray}

Formulas \eqref{52} and \eqref{54} demonstrate that the limits $\widehat \mu$ 
in \eqref{50} and $\widehat v$ in \eqref{55} satisfy the integral identity 
\eqref{40}, thus, in view of \eqref{57}, $\widehat v\in \cH (\eta)$ is a non-
trivial eigenfunction of the limit problem \eqref{40}, corresponding to its 
eigenvalue $\widehat \mu$.

\begin{lemma} \label{L31}
The limits $\widehat \mu$ in \eqref{50} and $\widehat v$ in \eqref{55} compose
an eigenpair of the limit problem \eqref{40}, or, \eqref{25}--\eqref{26}, \eqref{38}.
\end{lemma}

The proof of Lemma \ref{L31} is completed only in Remark \ref{remNEW}, but
on the other hand, it is needed only at the  of the proof of Theorem \ref{3AS}.

\subsection{Operator formulation of the model problem in the periodicity cell}
\label{sec3.2}
The next assertion is known as the "lemma on almost eigenvalues and eigenvectors" \cite{ViLu} and it is a consequence of the
spectral decomposition of the resolvent, see, e.g. \cite[Ch.\,6]{BiSo}.

\begin{lemma}\label{L32}
Let the function $\cU^h \in \cH_\eta^h$ and number $k^h \in \bbR_+$ be
such that 
\begin{eqnarray}
\Vert \cU^h ; \cH_\eta^h\Vert  = 1 \ \ \mbox{and} \ \ 
\Vert \cK_\eta^h \cU^h - k^h \cU^h ; \cH_\eta^h\Vert  =: \delta^h \in (0,
k^h). \label{62}
\end{eqnarray}
Then, there exists an eigenvalue $\kappa_m^h$ of the operator
$\cK_\eta^h$ such that 
\begin{eqnarray}
|k^h - \kappa_m^h| \leq \delta^h. \label{63}
\end{eqnarray}
Furthermore, for any $\delta_*^h  \in (\delta^h, k^h)$ one can find
a column of coefficients $b^h = \big( b_{M^h}^h, \ldots,  b_{M^h + X^h -1}^h
\big)$ satisfying
\begin{eqnarray}
\Big\Vert \cU^h - \sum_{p= M^h}^{M^h +X^h-1} b_p^h u_p^h;
\cH_\eta^h \Big\Vert  \leq \frac{\delta^h}{\delta_* ^h} , 
\ \ \ \sum_{p= M^h}^{M^h +X^h-1} | b_p^h|^2 = 1  , \label{64}
\end{eqnarray}
where the sequence $\kappa_{M^h}^h , \ldots , \kappa_{M^h+ X^h -1}^h$ 
consists of  all eigenvalues of $\cK_\eta^h$ contained in the closed 
interval $[k^h -\delta_*^h, k^h + \delta_*^h]$ and $u_{M^h}^h , \ldots , 
u_{M^h+ X^h -1}^h$
are the corresponding eigenvectors, which are orthonormalized in $\cH_\eta^h$.
\end{lemma}

\subsection{Global approximation and estimation of the discrepancies}
\label{sec3.3}
To glue the outer and inner expansion constructed in Section \ref{sec2}
we employ a trick, which uses cut-off functions with overlapping 
supports, see \cite[Ch.\,2]{MaNaPl} and \cite{na239} and Fig.\,\ref{fig2},b). 
Consequently, we define the following smooth functions
\begin{eqnarray}
& & X^h(x) = 1  \ \ \mbox{for } x \in \varpi_\sharp^h, \ n < - 2h , 
\row
X^h(x) =  0  \ \ \mbox{for } x \in \Theta \ \mbox{or} \ 
x \in \varpi_\sharp^h, \ n > - h ; \label{65}   \\
& & \chi(x) =  1  \ \ \mbox{for } x \in \varpi_\sharp^h \cap \cV, 
\ n >  - r_0 \ \ \mbox{or} \ x \in \Theta \cap \cV, \ r < r_ 0, 
\row
\chi(x) =  0  \ \ \mbox{for } x \in \varpi_\sharp^h , 
\ n <  - 2 r_0 \ \ \mbox{or} \ x \in \Theta, \ r > 2 r_ 0, 
\label{66}   \\
& & \cX^h(x) = 1  \ \ \mbox{for } x \in \Theta, \ r> 2 H_\theta h , 
\row
\cX^h(x) =  0  \ \ \mbox{for } x \in \varpi_\sharp^h  \ \mbox{or} \ 
x \in \Theta, \ r < H_\theta  h . \label{67}   
\end{eqnarray}
Here, $r_0$ is a positive number so small that $\{ x \in \varpi^h \, : 
\, r < 2 r_0 \} \subset \cV$  and $H_\theta = \max\limits_{ y \in \partial
\theta } H(y)$ and will be finally fixed only in Lemma \ref{L39}. 

If $\{ \mu(\eta), v(\eta; \cdot)\}$ is an eigenpair of the limit problem
\eqref{25}--\eqref{26}, \eqref{38},  formulas  \eqref{60} and \eqref{21} 
suggest to  take 
\begin{eqnarray}
k^h(\eta) = h^{-1} \big( \mu (\eta) +1 \big)^{-1}   \label{68}
\end{eqnarray}
as the "almost eigenvalue" of the operator $\cK_\eta^h$. In the following we
will omit $\eta$ from the notation. 

We take the following composite function as an "almost eigenvector":
\begin{eqnarray}
\cW^h &= & X^h(x) \big( v^0(y) + h^2 v'(y, h^{-1} z) \big)
\rowpl
\chi(x) \Big( a^0 + hH(y) \partial_n v^0(s,0) 
W\big(h^{-1} H(y)^{-1}n, h^{-1} H(y)^{-1} z  \big)- ha_h^1(s) \Big)
\rowmi
X^h(x)\chi(x)\big(v^0(s,0) + n\partial_n v^0(s,0)\big), \ \ 
x\in \varpi_\sharp^h ; 
\label{69}  
\\
\cW^h &= & \cX^h(x) \big( V^0 + h V'(x) \big)
\rowpl
\chi(x) \Big( a^0 + hH(s) \partial_n v^0(s,0) 
W\big( h^{-1} H(s)^{-1}n, h^{-1} H(s)^{-1} z  \big)- ha_h^1(s) \Big)
\rowmi
\cX^h(x)\chi(x)\big(a^0 + hH(s) \partial_n v^0(s,0) 
\frac2\pi \ln \frac{r}{h H(s)} \big), \ \ 
x\in \Theta . \label{70}
\end{eqnarray}
Let us explain the complicated structure of this global asymptotic
approximation of the eigenfunction $v(\eta, \cdot)$.

$\bullet$ In both expressions \eqref{69} and \eqref{70} the first two 
lines  on the right-hand sides contain terms which have been matched in 
Section \ref{sec2}, but subtracting the terms on the  third lines (with
factors $X^h \chi$ and $\cX^h h$, respectively), compensates such a 
duplication.

$\bullet$ The function $v'$ is found as a solution of the problem 
\eqref{24} and it is made unique by imposing the orthogonality condition
\begin{eqnarray*}
\int\limits_{- H(y)}^0 v'(y,\zeta) d\zeta = 0, \ \ y \in \omega_\sharp^0.
\end{eqnarray*}

$\bullet$ The arguments of the inner term \eqref{33} are different in
\eqref{69} and \eqref{70}, but nevertheless the function $\cW^h$
belongs to $H^1(\varpi^h)$ because $H(y)$ becomes $H(s)$ on the
common boundary surface \eqref{03} of $\varpi_\sharp^h$ and $\Theta$.

$\bullet$ The constant $a^0 = V^0$ is taken from \eqref{30} and \eqref{36}.
The behavior of $V'$ near the edge $\partial \theta$ is described
by the formula
\begin{eqnarray}
V'(x) = \frac2\pi H(s) \partial_n v(s,0) \ln r+ V_0'(s) + 
\widetilde V'(x) , \ \ \ x \in \Theta \cap \cV, \label{72N}
\end{eqnarray}
cf.\,\eqref{37}. Here, $V_0' \in C^\infty(\partial \theta)$ and 
\begin{eqnarray}
|\partial_z^k \partial_n^\ell \partial_s^m \widetilde V'(x)A \leq 
c_{k,\ell,m} r^{1-k-\ell} (1 + |\ln r|), \ \ x \in \Theta \cap \cV ; \label{72}
\end{eqnarray}
the orthogonality condition 
\begin{eqnarray*}
\int\limits_\Theta V'(x) dx = 0  
\end{eqnarray*}
makes both $V'$ and $V_0'$  uniquely defined. 

$\bullet$ The function $a_h^1$ in \eqref{70} is determined by
\begin{eqnarray}
a_h^1(s) = V_0'(s) - H(s) \partial_n v(s,0)
\frac{2 }{\pi}	\ln (hH(s)) . \label{73}
\end{eqnarray}

Recalling our definitions \eqref{69} and \eqref{65}, \eqref{66}, we see
that the function $\cW^h$coincides with $v(y) + h^2 v'(y,\zeta)$ near that
ends $\{ x \in \partial \varpi^h \, : \, y_1= \pm 1/2 \}$ of the cell
\eqref{01}, and therefore, satisfies the quasi-periodicity condition
\eqref{10} because \eqref{261}, \eqref{26} hold for $v$ as well as $v'$. 
Thus, the function $\cW^h$ lies in $H_\eta^1(\varpi^h)$ and we can normalize
it as
\begin{eqnarray}
\cU^h = \Vert \cW^h ; \cH_\eta^h \Vert^{-1} \cW^h. 
\label{73v}
\end{eqnarray}

Let us evaluate the quantity $\delta^h$, \eqref{62}, defined by a couple 
$\{ k^h, \cU^h\}$, using \eqref{58}, \eqref{59} and \eqref{68}
\begin{eqnarray}
\ \ \ \ \delta^h &=& \Vert \cK_\eta^h \cU^h - k^h \cU^h ; H_\eta^h \Vert 
=  \sup \big|\langle \cK_\eta^h \cU^h - k^h \cU^h, \cV^h 
\rangle_h \big| 
\roweq
k^h \Vert \cW^h ; \cH_\eta^h \Vert^{-1} \sup \big| (h + h\mu) 
( \cW^h, \cV^h)_{\omega^0} 
\rowmi
(\nabla \cW^h , \nabla \cV^h)_{\varpi^h} - h(\cW^h, \cV^h)_{\omega^0}
\big| 
\roweq
h^{-1} (1 + \mu)^{-1}  \Vert \cW^h ; \cH_\eta^h \Vert^{-1} \sup \big| 
(\nabla \cW^h , \nabla \cV^h)_{\varpi^h} 
\rowmi 
h \mu (\cW^h, \cV^h)_{\omega^0}
\big| , \label{74}
\end{eqnarray}
where the supremum is computed over the unit ball of $\cH_\eta^h$ so that 
we have
\begin{eqnarray}
\Vert \cV^h ; \cH_\eta^h \Vert \leq 1 . \label{75}
\end{eqnarray}
Our aim is to show that $\delta^h \leq c ( 1 + |\ln h|)$. 
This will follow in Corollary \ref{corNEW} from Lemmas \ref{L30} and \ref{L39}.

\begin{lemma} 
\label{L30} We have 
\begin{eqnarray}
\big| 
(\nabla \cW^h , \nabla \cV^h)_{\varpi^h} - h \mu (\cW^h, \cV^h)_{\omega^0}
\big| \leq c h^{3/2} ( 1 + |\ln h|)   \label{75a}
\end{eqnarray}
for all $\cV^h$ belonging to the unit ball of $\cH_\eta^h$. 
\end{lemma}

Proof. To evaluate the right-hand side of \eqref{74}, we write
\begin{eqnarray}
& & ( \nabla \cW^h , \nabla \cV^h)_{\varpi^h} - h \mu (\cW^h, \cV^h)_{\omega^0} 
\roweq
I_\varpi(\cV^h) + I_\Theta (\cV^h)+ I_\Xi(\cV^h)+ I_a(\cV^h) , \label{77}
\end{eqnarray}
where the terms will be defined in the following. First of all, we set
\begin{eqnarray}
I_\varpi(\cV^h)  &=& \big( \nabla( v^0 + h^2 v'), \nabla (X^h\cV^h) 
\big)_{\varpi_\sharp^h }-h \mu( v^0 + h^2 v', X^h\cV^h )_{\omega_\sharp^h }
\rowpl 
\big( ( v^0 + h^2 v') \nabla X^h , \nabla \cV^h \big)_{\varpi_\sharp^h }
- \big( \nabla( v^0 + h^2 v'),\cV^h \nabla X^h \big)_{\varpi_\sharp^h }
\nonumber \\
&=:& I_\varpi^1(\cV^h) + I_\varpi^2(\cV^h) + I_\varpi^3(\cV^h) 
+I_\varpi^4(\cV^h)  . \label{78}
\end{eqnarray}
Here, we commuted the gradient operator $\nabla$ with the cut-off function
$X^h$. Observing that, by \eqref{65}, $X^h \cV^h = 0$ on $\upsilon^h$ and
integrating by parts yield
\begin{eqnarray*}
& & I_\varpi^1(\cV^h) + I_\varpi^2(\cV^h)
\roweq
- \big( \Delta ( v^0 + h^2 v'), X^h \cV^h\big)_{\varpi_\sharp^h}
+ \big(\partial_z - h \mu) ( v^0 + h^2 v'), X^h \cV^h\big)_{\omega_\sharp^0}
\rowpl 
\big( \partial_\nu ( v^0 + h^2 v'), X^h \cV^h\big)_{\varpi_\sharp^h}.
\end{eqnarray*}
The consideration in Section \ref{sec2.1} gives the relations
\begin{eqnarray*}
& & \Delta(v^0 + h^2 v') = h^2 \Delta_y v' , \ \ (\partial_z - h \mu)
(v^0 + h^2 v') = - h^3 \mu v',
\row
(1 + h^2 |\nabla_yH|^2)^{1/2} ( \partial_\nu v^0 + h^2 v') 
= -h^3 \nabla_y H \cdot \nabla_y v'. 
\end{eqnarray*}
Hence, using the evident estimate
\begin{eqnarray*}
& & \Vert \cV^h; L^2(\varpi_\sharp^h ) \Vert^2 + h 
\Vert \cV^h; L^2(\omega_\sharp^h ) \Vert^2 
\rowleq
c\big( h^2 \Vert \partial_z \cV^h; L^2(\varpi_\sharp^h ) \Vert^2 
+ h \Vert \cV^h; L^2(\omega_\sharp^0) \Vert^2 \big) \leq c ,
\end{eqnarray*}
see \eqref{75}, gives us
\begin{eqnarray}
& & |I_\varpi^1(\cV^h) + I_\varpi^2 (\cV)|
\rowleq c\Big( h^2 ({\rm mes}_3(\varpi_\sharp^h))^{1/2} 
\Vert \cV^h; L^2(\varpi_\sharp^h ) \Vert
+ h^3 ({\rm mes}_2(\omega_\sharp^0))^{1/2} 
\Vert \cV^h; L^2(\omega_\sharp^0 ) \Vert
\rowpl
h^2 ({\rm mes}_2(\omega_\sharp^h))^{1/2} 
\Vert \cV^h; L^2(\omega_\sharp^h ) \Vert \Big)
\leq c \big( h^2 h^{1/2} + h^3h^{-1/2} + h^3 h^{-1/2} \big)
\roweq 
3c h^{5/2}. \label{790}
\end{eqnarray}
The remaining terms $I_\varpi^3(\cV^h)$ and  $I_\varpi^4(\cV^h)$
will be taken into account later. 

On the massive domain $\Theta$ we set 
\begin{eqnarray*}
\cV^h(x) = \cV_0^h + \cV_\perp^h(x), \ \ \mbox{where} \ \int\limits_\Theta
\cV_\perp^h (x) dx = 0,  
\end{eqnarray*}
and using definition \eqref{58} together with the Poincar\'e and trace
inequalities, cf.\,\eqref{521}, we obtain
\begin{eqnarray}
h^{1/2} |\cV_0^h| + \Vert \cV_\perp^h ; H^1(\Theta) \Vert \leq c \Vert \cV^h; 
\cH_\eta^h\Vert \leq c. \label{792}
\end{eqnarray}
Then, we define and estimate
\begin{eqnarray}
I_\Theta(\cV^h)  &=& \big( \nabla( \cX^h ( V^0 + h V')), 
\nabla \cV^h \big)_{\Theta}-h \mu( \cX^h( V^0 + h V'), \cV^h )_{\theta}
\roweq
- \big( \Delta ( V^0 + h V'), \cX^h \cV_\perp^h\big)_{\Theta}
+ \big( \partial_\nu ( V^0 + h V'), \cX^h \cV_\perp^h\big)_{
\partial \Theta \setminus \theta}
\rowpl
\big( -\partial_\nu - h \mu) ( V^0 + h V'), \cX^h \cV_\perp^h\big)_{
\theta} - h \mu \big(  \cX^h( V^0 + h V'), \cV_0^h\big)_{
\theta}
\rowpl 
\big(  ( V^0 + h V') \nabla \cX^h, \nabla \cV_\perp^h \big)_{\Theta}
- \big( \nabla ( V^0 + h V'), \cV^h_\perp \nabla \cX^h \big)_\Theta
\nonumber \\
&=:& \sum_{j=1}^8 I_\Theta^j (\cV^h) . \label{80}
\end{eqnarray}
Recall that $V^0$ is constant and that, by  Section \ref{sec2},
\begin{eqnarray*}
& & \Delta V'= 0 \ \mbox{in} \ \Theta, \  \ \ \partial_\nu V' = 0 \ \mbox{in} 
\ \partial \Theta \setminus \theta, \ \ \ ( \partial_z - h \mu)(V^0 + hV') = 
- h^2 \mu V' \ \mbox{in} \ \theta,  \row
|V'(y,0) | \leq c(1 + |\ln r|) \leq c ( 1 + |\ln h|) \ \mbox{for} \
y \in \theta \cap {\rm supp}\, \cX^h, 
\end{eqnarray*}
hence, we have 
\begin{eqnarray*}
& & I_\Theta^1(\cV^h) = 0, \ \ I_\Theta^2 (\cV^h) = 0, \row
|I_\Theta^3 (\cV^h) | \leq c h^2(1 + |\ln h|) \Vert \cV_\perp^h ; L^2(\theta) 
\Vert \leq c h^2 ( 1 + |\ln h |) . 
\end{eqnarray*}
The last two terms will be treated in the next paragraph, while the estimate
\eqref{792} yields
\begin{eqnarray*}
& & \big| I_\Theta^4(\cV^h) - h \mu V^0 \overline{\cV_0^h} 
{\rm mes}_ 2(\theta) \big| 
\rowleq c h \big({\rm mes}_ 2 (\theta_\cX^h) |V^0|
+ h \Vert V' ; L^2(\theta) \Vert \big) |\cV_0^h| 
\rowleq
ch(h+h) h^{-1/2} \Vert \cV^h ; \cH_\eta^h \Vert \leq 2 c h^{3/2},
\end{eqnarray*}
where, according to \eqref{67},
\begin{eqnarray*}
\theta_\cX^h = \{ y \in \theta \, : \, \cX^h(y,0) \not= 1 \} 
\subset \{ y \in \theta \, : \, r < 2 H_\theta h \}. 
\end{eqnarray*}

The last couples of terms in \eqref{78} and \eqref{80} will be 
compensated by parts of the expression
\begin{eqnarray}
& & I_a ( \cV^h) = 
\big( \nabla (X^h \chi ( v_\theta^0 + n \partial_z v_\theta^0) ),
\nabla v^h \big)_{\varpi_\sharp^h} - h \mu 
\big( X^h \chi ( v_\theta^0 + n \partial_n v_\theta^0) ),
\cV^h \big)_{\omega_\sharp^0}
\rowpl
\big( \nabla(\cX^h \chi(a^0 + hA_h-ha^1), \nabla \cV^h\big)_\Theta
\rowmi 
h \mu \big( \cX^h \chi(a^0 + hA_h-ha^1), \cV^h\big)_\theta, \label{821}
\end{eqnarray}
where $\partial_n v_\theta (s) = \partial_n v^0(s,0)$ and 
$A_h(s,r) = H(s) \partial_n v_\theta^0(s) \frac{2}{\pi} \ln \frac{r}{hH(s)}$.
Taking into account the position of the supports in \eqref{65}--\eqref{67}
shows that
\begin{eqnarray*}
\nabla( X^h \chi) = \nabla X^h + \nabla \chi,  \ \ \ 
\nabla( \cX^h \chi) = \nabla \cX^h + \nabla \chi,  
\end{eqnarray*}
and we rewrite \eqref{821} as follows:  
\begin{eqnarray}
& & I_a(\cV^h) \roweq
- \big( \nabla( v_\theta^0 + n \partial_n v_\theta^0) , 
\nabla( X^h \chi \cV^h) \big)_{\varpi_\sharp^h} 
+ h \mu \big(  v_\theta^0 + n \partial_n v_\theta^0 , 
X^h \chi \cV^h \big)_{\omega_\sharp^h} 
\rowmi
\big( ( v_\theta^0 + n \partial_n v_\theta^0) \nabla X^h , \nabla \cV^h 
\big)_{\varpi_\sharp^h} +
\big( \nabla( v_\theta^0 + n \partial_n v_\theta^0) , 
\cV^h \nabla X^h  \big)_{\varpi_\sharp^h}
\rowmi
\big( ( v_\theta^0 + n \partial_n v_\theta^0) \nabla \chi, \nabla \cV^h 
\big)_{\varpi_\sharp^h} +
\big( \nabla( v_\theta^0 + n \partial_n v_\theta^0) , 
 \cV^h\nabla \chi \big)_{\varpi_\sharp^h}
\rowmi
\big(\nabla(a^0 + hA_h -ha_h^1 ), \nabla(\cX^h \chi \cV_\perp^h)\big)_\Theta 
+ h \mu \big(a^0 + hA_h -h a_h^1 , \cX^h \chi \cV^h \big)_\theta
\rowmi
\big( (a^0 + hA_h -ha_h^1 ) \nabla\cX^h ,\nabla \cV_\perp^h\big)_\Theta 
+  \big( \nabla( a^0 + hA_h -h a_h^1) , \cV^h  \nabla \cX_\perp^h  \big)_\Theta
\rowmi
\big(\nabla(a^0 + hA_h -ha_h^1 ) \nabla \chi , \nabla \cV_\perp^h\big)_\Theta 
+ \big(\nabla ( a^0 + hA_h -h a_h^1 )  ,  \cV^h \nabla \chi \big)_\Theta
\nonumber \\
&=:& \sum_{j=1}^{12} I_\Theta^j(\cV^h) .  \label{82}
\end{eqnarray}
In view of \eqref{65} and \eqref{67} we have 
\begin{eqnarray*}
& & n = O(h) \ \mbox{in} \ T^h = {\rm supp}\,|\nabla X^h|, \ \ 
r = O(h) \ \mbox{in} \ \cT^h = {\rm supp}\,|\nabla \cX^h|, 
\row
|\nabla X^h(x) | + |\nabla \cX^h(x)| \leq c h^{-1}, \ \ \ 
{\rm mes}_3 T^h + {\rm mes}_3 \cT^h \leq c h^2 . 
\end{eqnarray*}
Hence, using the relations
\begin{eqnarray*}
& & | v^0(y) - v_\theta^0(s) - n \partial_n v_\theta^0(s) |\leq 
c ({\rm dist}\, (y,\partial \theta))^2, \ \ |v'(y,\zeta)| \leq c, 
\row
V^0 = a^0, \  \ \big|V^0(x) + h V'(x) - a^0- hA_h(s,r) + ha_h^1(s) \big|
\leq c \, {\rm dist}\, (x ,\partial \theta), 
\end{eqnarray*}
we conclude that 
\begin{eqnarray*}
& & |I_\varpi^3( \cV^h) + I_a^3(\cV^h)| 
\leq ch^{-1} ( {\rm mes}_ 3 T^h)^{1/2} h^2 \Vert \nabla \cV^h ; 
L^2(T^h) \Vert \leq ch^2 , 
\row
|I_\Theta^5( \cV^h) + I_a^9(\cV^h)| 
\leq ch^{-1} ( {\rm mes}_ 3 \cT^h)^{1/2} h^2 \Vert \nabla \cV^h ; 
L^2(\cT^h) \Vert \leq ch^2 . 
\end{eqnarray*}
The estimates
\begin{eqnarray}
& & \ \ \ \ \ \ |I_\varpi^4( \cV^h) + I_a^4(\cV^h)| 
\leq ch^{-1} ( {\rm mes}_ 3 T^h)^{1/2} h \Vert  \cV^h ; 
L^2(T^h) \Vert \leq ch^{3/2}   
\nonumber  \\
& & \ \ \ \ \ \ 
|I_\Theta^6( \cV^h) + I_a^{10}(\cV^h)| 
\leq ch^{-1} ( {\rm mes}_ 3 \cT^h)^{1/2} h \Vert \cV_\perp^h ; 
L^2(\cT^h) \Vert \leq  
\rowleq 
c h^{3/2} (1 + |\ln h|) 
\label{H3}
\end{eqnarray}
are obtained in the same way, but with the help of the next lemma. 
The remaining terms in\eqref{82} will be considered later, but notice
that until now all terms of \eqref{78} and \eqref{80} have
been taken into account. 

\begin{lemma}
\label{L34}
The following inequalities hold true:
\begin{eqnarray}
& & \Vert \cV^h ; L^2(T^h) \Vert \leq c h^{1/2} \Vert \cV^h; \cH_\eta^h\Vert,
\row
\Vert \cV_\perp^h ; L^2(\cT^h) \Vert \leq c h^{1/2}
(1 + |\ln h|)  \Vert \cV_\perp^h; H^1(\Theta) \Vert. \label{H4}
\end{eqnarray}
\end{lemma}

Proof of Lemma \ref{L34}. First of all, repeating the calculation \eqref{04} with an evident
modification yields 
\begin{eqnarray}
\Vert \cV^h ; L^2(\varpi_\sharp^h) \Vert^2 &\leq &
c\big( h \Vert \cV^h ; L^2(\omega_\sharp^0) \Vert^2 
+ h^2 \Vert \partial_z \cV^h ; L^2(\varpi_\sharp^h) \Vert^2 
\leq \Vert \cV^h ; \cH_\eta^h \Vert^2 . \label{H5}
\end{eqnarray}
Then, we write the Newton-Leibnitz formula in the form
\begin{eqnarray}
\cV^h (s,n,z) = \int\limits_{n_0}^n \frac{\partial}{\partial t}
\big( \chi_\varpi (t) \cV^h(s,t,z) \big) dt ,         \label{H6}
\end{eqnarray}
where $n_0 > 0$ is fixed such that $\{ x \in \varpi_\sharp^h \, : \,
n > - n_0\} \subset \upsilon$ and $\chi_\varpi$ is a smooth cut-off
function such that
\begin{eqnarray*}
\chi_\varpi (n) = 1 \ \mbox{for} \ n > - n_0/3, \ \ \ 
\chi_\varpi (n) = 0 \ \mbox{for} \ n < - 2n_0/3 .
\end{eqnarray*}
Integrating \eqref{H6} over $(-2h,0) \ni n$, $\partial \theta \ni s
$ and $(- hH(y),0) \ni z$, in other words over $T^h \ni x$, leads
to the inequality
\begin{eqnarray*}
\int\limits_{T^h} |\cV^h (x) |^2 dx \leq ch \int\limits_{\varpi_\sharp^h}
\big( |\nabla \cV^h (x)|^2 + |\cV^h(x)|^2 \big) dx,
\end{eqnarray*}
which is nothing but the first inequality in \eqref{H4}.

The second inequality in \eqref{H4} follows from the estimate
\begin{eqnarray*}
\Vert \cV _\perp^h ; H^1(\Theta) \Vert \leq c \Vert \cV^h; \cH_\eta^h \Vert,
\end{eqnarray*}
cf.\,\eqref{521}, and the one-dimensional Hardy inequality with logarithm
\begin{eqnarray*}
\int\limits_0^{r_0} \Big| \ln  \frac{r}{r_0} \Big|^{-2} |\bfV (r) |^2 
\frac{d r}{r} \leq 4 
\int\limits_0^{r_0} \Big|  \frac{d \bfV}{dr} \Big|^{2} r dr 
\end{eqnarray*}
and the fact that $r^{-2} |\ln r|^{-2}  \geq ch^{-2} ( 1 + |\ln h|)^{-2}$ 
for $x \in \cT^h$. \ \ $\boxtimes$

\bigskip

We complete the proof of Lemma \ref{L30}. 
Let us write the term $I_\Xi ( \cV^h)$ as follows:
\begin{eqnarray}
I_\Xi ( \cV^h) &=& 
\big( \nabla \bfW_\varpi, \nabla(\chi \cV^h) \big)_{\varpi_\sharp^h}
- h \mu \big( \bfW_\varpi, \chi \cV^h \big)_{\omega_\sharp^0}
\rowpl
\big( \nabla W_\Theta, \nabla(\chi \cV^h) \big)_{\Theta}
- h \mu \big(  \bfW_\Theta, \chi \cV^h \big)_{\theta}
\rowpl
\big(  \bfW_\varpi \nabla \chi, \nabla \cV^h \big)_{\varpi_\sharp^h} 
- \big( \nabla \bfW_\varpi, \cV^h \nabla\chi  \big)_{\varpi_\sharp^h}
\rowpl
\big(  \bfW_\varpi \nabla \chi, \nabla \cV^h \big)_{\Theta} -
\big( \nabla \bfW_\Theta, \cV^h \nabla \chi  \big)_{\Theta}
= \sum_{k=1}^8 I_\Xi (\cV^h) , \label{H7}
\end{eqnarray}
where $\bfW_\varpi$ and $\bfW_\Theta$ are the multipliers of
$\chi$ in the middle of the right-hand sides in \eqref{69} and \eqref{70},
respectively; recall that the arguments of $W$ in these sums differ from 
each other. 

We treat the sums \eqref{H7} and \eqref{82} together and start with the
simplest terms. According to \eqref{33} and \eqref{31}, the modulus of the 
difference
\begin{eqnarray}
& & h H(y) \partial_n v^0(s,0) W \big( h^{-1} H(y)^{-1} n, h^{-1} H(y)^{-1} 
z \big) - n \partial_n v^0(s,0) 
\roweq
h H(y) \partial_n v^0(s,0) c_\Xi + O( e^{- \pi |n|/ (hH(y) } )
\label{H0}
\end{eqnarray}
is bounded by $ch$, we obtain
\begin{eqnarray}
\big| I_\Xi^5( \cV^h) + I_a^5(\cV^h)\big| 
\leq c h ({\rm mes}_3 \varpi_\sharp^h )^{1/2}  \Vert \cV^h ; 
L^2(\varpi_\sharp^h) \Vert  \leq c h^{3/2} . \label{H8}
\end{eqnarray}
In a similar way, by \eqref{72N} and \eqref{72} we observe that
the expression
\begin{eqnarray}
& & h H(s)\partial_n  v^0(s,0) W\big( h^{-1} H(s)^{-1} n ,  
h^{-1} H(s)^{-1} z \big)
\rowmi
ha_h^1(s) - h H(s) v^0(s,0) \frac{2 }{\pi	} \ln \frac{r}{hH(s)}
\label{H9}
\end{eqnarray}
is $O\big(h^2 r (1 + |\ln r|) \big)$ in $\Theta$ and thus $O(h^2)$ in 
supp\,$|\nabla \chi|$. Hence,
\begin{eqnarray*}
\big| I_\Xi^7( \cV^h) + I_a^{11}(\cV^h)\big| 
\leq c h^2  \Vert \cV_\perp^h ; L^2(\Theta) \Vert  \leq c h^2 . 
\end{eqnarray*}
As for the pairs $I_\Xi^6( \cV^h)$, $I_a^6( \cV^h)$ and 
$I_\Xi^8( \cV^h)$, $I_a^{12}( \cV^h)$, we have to take into account the
formula 
\begin{eqnarray}
\nabla = \big( (1 + n \varkappa(s))^{-1}\partial_s, \partial_n, \partial_z 
\big) .  \label{H11}
\end{eqnarray}
We observe that none of the three derivatives of \eqref{H11} affects the
above presented boundedness for the expression \eqref{H0}. The same holds concerning 
the estimate for the bounded of the expression \eqref{H9}, except on 
the set supp\,$|\nabla \chi|$, the position of which is however at
a positive distance from the edge $\partial \theta$ containing the 
singularities of the solutions. Thus, we conclude that 
\begin{eqnarray}
& & \big| I_\Xi^6( \cV^h)+ I_a^6( \cV^h)\big| \leq c h^{3/2} 
\Vert \cV^h ; L^2(\varpi_\sharp^h) \Vert \leq ch^{3/2} 
\row
\big| I_\Xi^8( \cV^h)+ I_a^{12}( \cV^h)\big| \leq c h^2 
\Vert \cV^h ; L^2(\Theta) \Vert \leq ch^{3/2}. \label{H12}
\end{eqnarray}
On the first line we also used \eqref{H5} and on the second one
the inequality
\begin{eqnarray*}
\Vert \cV^h; L^2(\Theta) \Vert \leq c\big( \Vert \nabla \cV^h; L^2(\Theta) 
\Vert + \Vert \cV^h; L^2(\theta) \Vert \big)   \leq
ch^{-1/2} \Vert \cV^h; \cH_\eta^h \Vert
\end{eqnarray*}
coming from the definition of the scalar product in $\cH_\eta^h$.

We are left with the sums
\begin{eqnarray}
I_\Xi^1 ( \cV^h) + I_\Xi^3( \cV^h), \ \ I_\Xi^2( \cV^h) + I_\Xi^4( \cV^h)
\label{H14}
\end{eqnarray}
and 
\begin{eqnarray}
I_a^1( \cV^h) + I_a^7( \cV^h), \ \ I_a^2( \cV^h) + I_a^8( \cV^h) 
\label{H15}
\end{eqnarray}
in \eqref{H7} and \eqref{82}, respectively. The two-dimensional integrals
included into the second  sums in \eqref{H14} and \eqref{H15}  are treated
using the information we already have on the expressions \eqref{H0} 
and \eqref{H8}, which gives us
\begin{eqnarray*}
& & \big| I_\Xi^2( \cV^h)+ I_a^2( \cV^h)\big| \leq c h^2
\Vert \cV^h ; L^2(\omega_\sharp^0) \Vert \leq ch^{3/2} 
\row
\big| I_\Xi^4( \cV^h)+ I_a^8( \cV^h)\big| \leq h^2 
\Vert \cV^h ; L^2(\theta) \Vert \leq ch^{3/2} . 
\end{eqnarray*}
Here we also applied the inequality $\Vert \cV^h; L^2(\omega^0) \Vert
\leq h^{-1/2} \Vert \cV^h ; \cH_\eta^h \Vert $, which is a consequence
of the definition \eqref{58}. 

To conclude with the first sums in  \eqref{H14} and \eqref{H15}, we
need to take into account several aspects. To start with, the first component 
on the gradient operator \eqref{H11}, which acts in the integrals over the
sets $\varpi_\sharp^h$ and $\Theta$, is treated in the same way as in
\eqref{H12}; notice that in the bound $ch^2 r(1+ |\ln r|)$ of the modulus
of \eqref{H9}, the factor $r(1+ |\ln r|)$ is small near the edge $\partial
\theta$. Next, the factor $\ n \varkappa(s)$ in the Jacobian of the
differential
\begin{eqnarray*}
dx = ( 1 + n \varkappa (s) ) ds dn dz 
\end{eqnarray*}
is also $O(r)$. Together, these compensate the singularities of 
$\nabla_\xi W(\xi)$ of the solution \eqref{33} of the Neumann problem 
\eqref{32} at the corner point $P$. 

Finally, recalling that the arguments of the function $W$ are different
in \eqref{69} and \eqref{70}, we write
\begin{eqnarray*}
& & \frac{\partial}{\partial z} \Big(
W \Big( \frac{n}{hH(y)} ,  \frac{z}{hH(y)} \Big) -  
W \Big( \frac{n}{hH(s)} ,   \frac{z}{hH(s)} \Big) \Big) 
= \frac1h   \Big( \frac{1}{H(y)} -  \frac{1}{H(s)} \Big)
\frac{\partial W }{ \partial \xi_2 } (\xi) ,
\row
\frac{\partial}{\partial n} \Big(
W \Big( \frac{n}{hH(y)} ,  \frac{z}{hH(y)} \Big) -  
W \Big( \frac{n}{hH(s)} ,   \frac{z}{hH(s)} \Big) \Big) 
\roweq
\frac1h   \Big( \frac{1}{H(y)} -  \frac{1}{H(s)} \Big)
\frac{\partial W }{ \partial \xi_1 } (\xi)
- \frac{n}h  \frac{\partial_s H(y) }{H(y)^2 } 
\frac{\partial W }{ \partial \xi_1 } (\xi)
- \frac{z}h  \frac{\partial_s H(y) }{H(y)^2 } 
\frac{\partial W }{ \partial \xi_2} (\xi) .
\end{eqnarray*}
Again, in $\varpi_\sharp^h$ there appears a factor of order $O(n)$,
due to $|H(y) - H(s)| = O(n)$, see Section \ref{sec2.2}. This and the
obvious relation
\begin{eqnarray*}
\int\limits_0^1 n^2 e^{-2 \delta n /h} dn = O(h) , \ \ \ \delta > 0,
\end{eqnarray*}
yield  an additional factor $h^{1/2}$ in the above estimates, thus
compensating the coefficient $h^{-1}$ in \eqref{31}.

We denote by $\widehat I_\Xi^q(\cV^h) $ and $\widehat I_a^q(\cV^h) $
the expressions in \eqref{H14} and \eqref{15} after the above mentioned
simplifications have been made. We have 
\begin{eqnarray}
& & \widehat I_\Xi^1(\cV^h) + \widehat I_\Xi^3 (\cV^h)
\roweq
\int\limits_{\partial \theta} \int\limits_\Xi \nabla_\xi
\big( a^0 + hH(s) \partial_n v^0(s,0) W(\xi)
-h a_h^1(s) \big) \cdot \nabla_\xi (\chi \cV^h) d \xi ds ,
\row
\widehat I_a^2(\cV^h)
= - \int\limits_{\partial \theta} \int\limits_\bbP  \nabla_\xi
\big( v^0(s,0) + h \xi_1 \partial_n v^0(s,0)  \big) \cdot 
\nabla_\xi (\chi X^h \cV^h) d \xi ds ,
\row
\widehat I_a^7(\cV^h) 
\roweq
- \int\limits_{\partial \theta} \int\limits_\bbK \nabla_\xi
\Big( a^0 + hH(s) \partial_n v^0(s,0) \frac{2}{\pi} \ln 
\frac{r}{hH(s) } \Big) \cdot \nabla_\xi (\chi  \cX^h \cV^h) d \xi ds .
\label{H16}
\end{eqnarray}
These integrals involve several solutions of the Neumann Laplacian
in the semi-strip $\bbP$, the quadrant $\bbK$ and in their union
$\Xi$, see Section \ref{sec2.2}, as well as test functions
with compact support, which in particular vanish near the end of 
$\bbP$ and the corner point of $\bbK$. As a consequence, the 
integrals \eqref{H16} vanish. 

This completes the consideration of the terms in \eqref{78}, \eqref{80},
\eqref{82} and \eqref{H7}.  Lemma \ref{L30} thus follows by  
combining the formulas \eqref{74}--\eqref{78} with the inequalities 
\eqref{790}--\eqref{H16}. \ \ $\boxtimes$

\bigskip

Note that the worst (largest) bound for the left-hand side of
\eqref{75a}  appears in \eqref{H3}. 

We also need the following estimates for the
approximate eigenfunctions.

\begin{lemma} \label{L39}
Let $\mu_m(\eta)$, $m \in \bbN$, be the eigenvalues \eqref{41} and 
let $v_m^0(\eta; \cdot)$  be the corresponding
eigenfunctions of the problem \eqref{40}, orthonormalized as in \eqref{43}.
Then, the functions $\cW_m^h$, which are defined in \eqref{69}, 
\eqref{70} by using $v_m^0$ and small enough $r_0 = r_0(m)> 0$, satisfy 
for some $\widetilde h_m > 0$ the relation
\begin{eqnarray}
\Vert   \cW_m^h - v_m^0 ; L^2(\omega^0) \Vert \leq 2^{-m-2}
\ \ \ \forall \, h \leq \widetilde h_m.    \label{9.0}
\end{eqnarray}
Consequently, for every $m \geq 2$ and $0 < h <  \widetilde h_m$, 
the sequence $(\cW_n^h)_{n=1}^m$  
is linearly independent, and there also holds, for some constants 
$c_1, c_2 > 0$, 
\begin{eqnarray}
c_1 h^{1/2} \leq \Vert \cW_m^h ; \cH_\eta^h \Vert \leq c_2 h^{1/2}
\ \ \ \forall \, h \leq \widetilde h_m.   \label{9.1}
\end{eqnarray} 
\end{lemma}

Proof.  
Let us consider the defining formulas \eqref{69}, \eqref{70}. 
By standard elliptic regularity results, \cite{ADN}, the smoothness
of the boundary $\partial \omega_\sharp^0$, and the definition as the solution of
\eqref{25}--\eqref{26}, \eqref{38}, the functions $v^0$, 
$\partial_n v^0$, $\Delta_y v^0$  are bounded by a constant independent of $h$ 
and $\eta$ in the domain $\omega^0$ (we again drop the index $m$ from the 
notation). Consequently, by \eqref{24},
also $v'$ is bounded in the same way. 

By \eqref{33}, we have $W(h^{-1} H(s)^{-1}n , h^{-1} H(s)^{-1} z) 
\leq c h^{-1}$, hence,
\begin{eqnarray*}
\big| h H(y) \partial_n v^0(s,0) W \big(h^{-1} H(s)^{-1}n , h^{-1} H(s)^{-1} z 
\big)  \big| \leq c'
\end{eqnarray*}
on $\omega^0$. The modulus of the function $a_h^1(s)$ is bounded by
$c h |\ln h|$, see \eqref{73}. Thus, the second and third rows of 
\eqref{69} can be written, by taking $ x = (y,0)$, as 
$ 
\chi(y) F^h(y)$, $ y \in \omega_\sharp^0,
$ 
where $F^h$ is a function, which is bounded by a constant independent of $h$. 
Choosing $r_0 = r_0(m)> 0$  small enough in the definition of $\chi$, \eqref{66},  we get
\begin{eqnarray*}
\Vert \chi F^h ; L^2 (\omega_\sharp^0) \Vert \leq 2^{-m-4}  .
\end{eqnarray*}
In the same way, the second and third rows of 
\eqref{70} can be written  as 
$ 
\chi(y) F^h(y)$, $ y \in \theta,$
for a  bounded extension of $F^h$ into $\theta$. Diminishing the number $r_0 > 
0$, if necessary,  yields also 
\begin{eqnarray}
\Vert \chi F^h ; L^2 (\theta) \Vert \leq  2^{-m-4}  \ \ \
\Rightarrow \ \ \ \Vert \chi F^h ; L^2 (\omega^0) \Vert \leq  2^{-m-3} 
 \label{76c}
\end{eqnarray}
for all $h$.  
Finally, in the subdomain $\widetilde \omega^h = \{ y \in \omega^0 \, :
\, X^h(y) = 1 \ \mbox{or} \ \cX^h(y) = 1 \}$ we have 
\begin{eqnarray}
\cW^h = v^0 + \chi F^h + \widetilde F^h,  \label{76b}
\end{eqnarray} 
where
$\widetilde F^h$ equals $h^2 v'(y, 0)$ on $\omega_\sharp^0$ and 
$h V'$ on $\theta$ so that $\Vert \widetilde F^h ; L^2 (\omega^0) \Vert 
\leq ch$. Since the area of $\omega^0 \setminus  \widetilde \omega^h$
is $O(h)$ and both functions $\cW^h$ and $v^0$ are bounded, the relation
\eqref{9.0} follows from \eqref{76c} and  \eqref{76b} by choosing 
$\widetilde h_m$ small enough.

The claim on the linear independence of the sequence follows from
the orthonormality relation \eqref{43} (note that the expression on the 
left-hand side of \eqref{43} is nothing but the natural inner product
of $L^2(\omega^0)$, since the eigenfunctions $v_m^0$ are
constant on $\theta$), \eqref{9.0} and a well-known result concerning 
perturbations of orthonormal bases.  

The lower bound in \eqref{9.1} follows 
from \eqref{43}, \eqref{9.0} and the definition of the norm of 
$\cH_\eta^h$ in \eqref{58}. The upper bound is a matter of
a straightforward calculation of the $\cH_\eta^h$-norm of the 
representation \eqref{69}--\eqref{70}. 
\ \ $\boxtimes$

\bigskip

We can now conclude with the estimation of $\delta^h$, \eqref{74}:
by \eqref{9.1},  the multiplier of the last $\sup | \ldots |$ in \eqref{74} 
does not exceed $c h^{-1} \Vert \cW^h ; \cH_\eta^h \Vert^{-1}
\leq C h^{-3/2}$. Thus, Lemma \ref{L30} yields

\begin{corollary}\label{corNEW} 
For $ 0< h < \widetilde h$, there holds the bound 
\begin{eqnarray}
\delta^h \leq h^{-3/2} h^{3/2} ( 1 + |\ln h|) = c ( 1 + |\ln h|) .
\label{H17}
\end{eqnarray}
\end{corollary}

Here, $ \widetilde h > 0$ means the number $ \widetilde h_m$, found in  
Lemma \ref{L39} for the approximate eigenvector under consideration.

\subsection{Theorem on asymptotics}  \label{sec3.4}

Let us state our main asymptotic result on the eigenvalues
of the model problem. 

\begin{theorem}
\label{3AS} For all $n \in \bbN$ there exist positive numbers 
$h_n$ and $c_n$ such that the entries of the eigenvalue sequences
\eqref{12} and \eqref{41} of the problems \eqref{7}--\eqref{11} and
\eqref{40}, respectively, are related  by
\begin{eqnarray}
|\Lambda_n^h(\eta) - h \mu_n(\eta) | \leq c_n h^2 (1+ |\ln h|)
\ \ \ \forall \, h\in (0, h_n], \ \eta \in [- \pi, \pi].   \label{vantaa}
\end{eqnarray}
\end{theorem}

As for the proof, one technical difficulty is caused by the possible
higher multiplicities of the eigenvalues of the limit problem. To treat
this we need to proceed in several steps and thus start by  showing 
the bound \eqref{vantaa} with $\mu_m $ in place 
$\mu_n$ for some $m \in \bbN$, which is for the moment unspecified.
In \eqref{63} of Lemma \ref{L32} and 
\eqref{74},  
\eqref{H17} we found an eigenvalue $\kappa_n^h = \kappa_n^h(\eta)$ of the model 
problem such that for an eigenvalue $k_m^h= k_m^h(\eta)$ of the limit problem  there
holds 
\begin{eqnarray*}
| \kappa_n^h  - k_m^h | \leq c(1 + |\ln h|)  .
\end{eqnarray*}
We have by \eqref{60} and \eqref{68}, 
\begin{eqnarray}
& & c(1 + |\ln h|) \geq |\kappa_n^h  - k_m^h  | = 
\Big| \frac{1}{\Lambda_n^h    + h } - \frac{1}{h \mu_m   + h } \Big|  
\end{eqnarray}
hence,
\begin{eqnarray}
& & 
 | \Lambda_n^{h}  -  h  \mu_m   | = | \Lambda_n^{h}  + h - 
 ( h  \mu_m   + h ) |
\rowleq 
c(1 + |\ln h|) \big( \Lambda_n^{h}  + h  \big)
\big(  h  \mu_m   + h \big) .  \label{vantaa2}
\end{eqnarray}
Here, both $\mu_m $ and $\Lambda_n^h$ are bounded by some
constant $c_{m,n} > 0$ so that \eqref{vantaa2}
and the triangle inequality give us
$$
|\Lambda_j^{h}  |\leq h \mu_m  +  c_{m,n}' (1 + |\ln h|) h 
\leq c_{m,n} (1 + |\ln h|) h .
$$
Inserting this again to the right hand side of 
\eqref{vantaa2} yields $|\Lambda_j^{h}  |\leq h \mu_m 
+ O(h^2 (1 + |\ln h|)^2) = O(h)$. 
Using this and \eqref{vantaa2} once more gives us
\begin{eqnarray}
|\Lambda_n^{h} - h  \mu_m |  
\leq  C_n h^2 (1 + |\ln h|),    \label{vantaa3}
\end{eqnarray}
i.e. \eqref{vantaa} holds true with $\mu_m $ in the place of $\mu_n$.

\begin{remark}
\label{remNEW}
We can now complete the proof for the formula \eqref{49} and thus
also Lemma \ref{L31}. Indeed, fixing $m \in \bbN$, we have above found for each  
eigenvalue $\mu_k$, $k \leq m$, with multiplicity $\varkappa_k$ at least
$\varkappa_k$ eigenvalues $\Lambda_{J(k)}^h, \ldots , 
\Lambda_{J(k) + \varkappa_k -1}$ belonging to the the interval
\begin{eqnarray}
\big[  h \mu_k - C_k h^2 (1 + | \ln h| ) , 
 h \mu_k + C_k h^2 (1 + | \ln h| ) \big]   \label{intX}
\end{eqnarray}
for $h \in (0, h^{(k)} )$. At this point we do not know, if there are 
eigenvalues outside the segments  \eqref{intX}, nevertheless, we can conclude
that, for $h < \min \{ h^{(1)}, \ldots, h^{(m)} \}$, there are at least 
$m-1$ eigenvalues in the intervals \eqref{intX} with $k = 1 , \ldots , m-1$.
Thus, $J(m) \geq m$ and, hence, 
$$
\Lambda_{m+j}^h(\eta) \leq \Lambda_{J(m)+j}^h(\eta) \leq 
 h \mu_m +  C_m h^2 (1 + |\ln h|)  \leq C_m' h ,\ \ \ j=0,\ldots, 
 \varkappa_m-1.
$$
%
\end{remark}

\begin{remark}
\label{remNEW1}
Let us  fix $n \in \bbN$ and $\eta \in [0, 2 \pi)$ and thus also 
the eigenvalue $\mu_n= \mu_n(\eta) $, see \eqref{41}, 
and denote its multiplicity by  $\cJ=\cJ(n,h) \in \bbN$. We claim that 
for some $\hat h > 0$, the interval 
\begin{eqnarray}
\vartriangle_{n,\eta} = \big[ h(\mu_n  - C_n ), h( \mu_n  + C_n) 
\big],
\label{9.2}
\end{eqnarray} 
contains for all $h \in (0, \hat h]$ at least  $\cJ$ eigenvalues 
$\Lambda_p^h(\eta)$ 
with multiplicities counted. In \eqref{9.2}, the number $C_n = C_n(\eta)$
is, say, half of the distance of $\mu_n(\eta)$ to the nearest
different eigenvalue of the limit problem.

Given $\cJ$ orthonormalized  eigenvectors $v_{p(n) + j}^0$,  $j = 1, \ldots 
, \cJ $, see \eqref{43}, 
corresponding to $\mu_n$  we construct $\cJ$ approximate
eigenvectors $\cW_{p(n) + j}^h$ for all  $j$ as in 
\eqref{69}, \eqref{70} by using the vectors $v_{p(n) + j}^0$. 
Lemma \ref{L39} then shows that the vectors $\cW_{p(j)}^h$, $j = 1, 
\ldots , \cJ $ form a linearly independent sequence, if $h $ is small enough.

For any $j$, we take  $\cU_{p(n) + j}^h = \cW_{p(n) + j}^h
\Vert \cW_{p(n) + j}^h ; H_\eta^h \Vert^{-1}$, see \eqref{73v},
and  $\delta_*^h= h^{-3/4}$ in Lemma  \ref{L32}; note that the latter 
choice is possible, 
since $\delta_*^h \in (\delta^h, k^h)$ and $k^h $ is of order $h^{-1}$. 
Due to \eqref{74} and the estimate \eqref{H17} in Corollary
\ref{corNEW},  
the inequality \eqref{64}   holds with 
$\delta^h / \delta_*^h = c h^{3/4} (1 + |\ln h |)$ 
and with the  true eigenvectors $U_{k(n) + j}^h$,
$j= 1, \ldots ,X^h$, see \eqref{16},  of $\cK_\eta^h$ 
in place of $u_p^h$, with the corresponding
eigenvalues $\kappa_{k(n) + j}^h$, \eqref{61}, belonging to the interval
$
[ k_n - \delta_*^h, k_n + \delta_*^h] .
$
We denote for a moment by $\Lambda$  and $\kappa$  any of the 
eigenvalues $\Lambda_{k(n) + j}^h$ and $\kappa_{k(n) + j}^h$, $j= 1, 
\ldots ,X^h$, respectively, and obtain by the relations \eqref{60}, \eqref{68} 
\begin{eqnarray*}
h^{-3/4} = \delta_*^h \geq | \kappa^h - k_n^h | = \Big| \frac{1 }{ \Lambda^h + h } - \frac{1}{h \mu_n + h}
\Big| = \frac{|\Lambda^h - h \mu|}{h(\Lambda^h + h )( \mu + 1)} .
\end{eqnarray*}
Hence, since $\Lambda^h \leq c h$ by  \eqref{49}, we get
\begin{eqnarray}
|\Lambda^h - h \mu_n| \leq ch^{5/4} 
\end{eqnarray}
and this shows that the eigenvalues 
$\Lambda_{k(n) + j}^h$, $j= 1, \ldots ,X^h$, 
belong to   $\vartriangle_{n,\eta}$. According to \eqref{58}, \eqref{64},
we get for some numbers $b_j$,
\begin{eqnarray}
& & \Big\Vert \cU_{p(n) + j}^h - \sum_{j = 1}^{X^h} b_j 
U_{p(n) + j}^h ; L^2(\omega^0) \Big\Vert 
\rowleq h^{-1/2}
\Big\Vert \cU_{p(n) + j}^h - \sum_{j = 1}^{X^h} b_j U_{p(n) + j}^h
; \cH_\eta^h \Big\Vert \leq h^{-1/2} \frac{\delta^h}{\delta_*^h}
\leq  ch^{1/5} .  \label{9.3}
\end{eqnarray}
Now, if $X^h < \cJ$, we arrive at a contradiction using elementary Hilbert
space geometry, since
given any at most $\cJ -1$-dimensional subspace $\cY$ of $L^2(\omega^0)$, 
one of the $\cJ$ nearly orthonormal  vectors $\cW_{p(n)+j}^h$, see  
\eqref{9.0}, \eqref{43}, is at least at the distance $ \frac12 
\cJ^{-1/2}$ of (any vector of) $\cY$. The same then holds for the 
corresponding $\cU_{p(n)+j}^h$, which contradicts \eqref{9.3},
where the span of the $X^h$ functions $U_{p(n) + j}^h$ is taken for $\cY$ .
\end{remark}

Proof of Theorem \ref{3AS}.
From Remark \ref{remNEW1} we actually obtain that given $n \in \bbN$, 
\eqref{vantaa3} holds for some $m \leq n$: note that the entries in both 
sequences  \eqref{12}  and \eqref{41} are in an increasing order and that
the intervals \eqref{9.2} are disjoint. 

It  thus 
suffices to show that $m \geq n$ in \eqref{vantaa3}, and to this
end it is enough to prove that  in the notation of Remark 
\ref{remNEW1}, the case $\cJ < X^h$ cannot happen; here, $\cJ$ is the multiplicity 
of $\mu_n$ and $X^h$ is the total multiplicity of the eigenvalues 
$\Lambda_{p(j)}^h$ belonging to  the interval $\vartriangle_{n,\eta}$, \eqref{9.2}.

Suppose  that for some $\widehat h > 0$ we have 
$\cJ < X^h$, for infinitely many $h \in (0, \widehat h)$ forming 
a set with 0 as an accumulation point. Then, for each such $h$
have $X^h$ eigenfunctions $U_{p(j)}^h  $,
$j = 1 , \ldots , \cJ$ which are orthogonal to each other in the 
norm of $L^2(\omega^0)$, see  \eqref{16}. According to 
Lemma  \ref{L31} and its proof,  and the choice of the interval 
$\vartriangle_{n,\eta}$  we find a  sequence $\{ h_q \}_{q=1}^\infty$ such that 
$$
h_q^{-1} \Lambda_{p(j)}^{h_q}(\eta) \to  \mu_n(\eta)  \ \ 
\mbox{as} \ q \to + 
\infty
$$
for all $j = 1, \ldots , X^h$, and also the eigenfunctions 
$U_{p(j)}^h $ converge to some eigenfunctions 
$v_{p(j) }^0  $ of the eigenvalue $ \mu_n(\eta) $ in the 
norm of $L^2(\omega_\sharp^0)$, as $q \to + \infty$; see \eqref{55}. 
This leads to the conclusion that the $X^h$ limit functions 
$v_{p(j) }^0 $ are also mutually orthogonal in 
$L^2(\omega_\sharp^0)$ so that they in particular are linearly
independent. Thus, there exist at least  $\cJ +1 $,
linearly independent eigenfunctions of the limit problem 
corresponding to $\mu_n $, which contradicts 
the choice of $\cJ$.

Now, the eigenvalues of both the model and limit problems are
arranged into an increasing order, so the indices $m$ and $n$
must be the same in \eqref{vantaa3}, i.e., \eqref{vantaa} holds. 
\ \ $\Box$

\section{Detecting spectral gaps} \label{sec4}
\subsection{Limit problem with a new small parameter} 
\begin{figure}
\begin{center}
\includegraphics[ height=6cm,width=15cm]{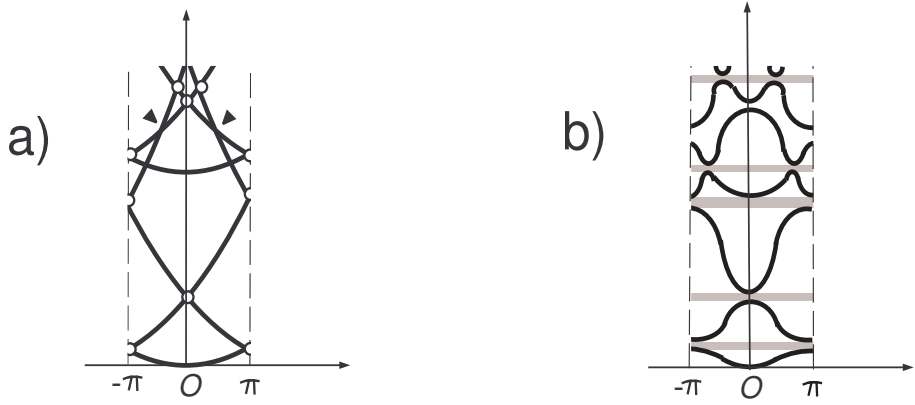}
\end{center}
\caption{a), b) Dispersion curves, c) opening of narrow gaps.}
\label{fig4}
\end{figure}

\label{sec{4.1}}
In order to analyse  the possible existence of spectral gaps in the spectrum \eqref{16}
of the Steklov problem \eqref{u1}--\eqref{u3}, we consider the limit
problem \eqref{25},\eqref{260},\eqref{26},\eqref{38} with the additional
assumptions that
\begin{eqnarray}
H(y) = 1 , \ \ \ \theta^\varepsilon = 
\{ y \, : \,  \zeta := \varepsilon^{-1} y 
\in \theta^1 \},  \label{X1}
\end{eqnarray}
where $\theta^1$ is a domain in $\bbR^2$ surrounded by a smooth 
simple closed contour $\partial \theta^1$ and $\varepsilon > 0$ 
is a new small parameter. We denote by $\{ v^\varepsilon(y ) , 
\mu^\varepsilon(\eta)\}$ an eigenpair of the limit problem in this
special case. Sending $\varepsilon \to +0$, which means glueing the small
hole, we end up with the following problem in the rectangle
$\omega = (-1/2, 1/2) \times (-\ell, \ell)$:
\begin{eqnarray}
& &  - \Delta_y v^0(\eta; y) 0  \mu^0(\eta) v^0(\eta; y), \,  y \in 
\omega, 
\row
\pm \frac{\partial v^0}{\partial y_2} (\eta; y_1, \pm \ell ) 
= 0, \ |y_1| < \frac12, 
\row
v^0\Big(\eta; \frac12, y_2   \Big)  =  e^{i \eta}
v^0\Big(\eta; - \frac{1}2, y_2  \Big),  \  
|y_2| < \ell \row
\frac{\partial v^0}{\partial y_1}\Big(\eta;\frac12, y_2  \Big) 
= e^{i \eta} \frac{\partial v^0}{\partial y_1}\Big(\eta;-\frac12, y_2  \Big),
\  |y_2| < \ell, 
 \label{X2} 
\end{eqnarray}
the solutions of which are given by the formulas
\begin{eqnarray*}
& & \mu_{jk}^0(\eta) = ( 2 \pi j + \eta) ^2 + \frac{\pi^2 k^2}{4 \ell^2 },
\row
v_{jk}^0(\eta; y) = e^{i(2 \pi j + \eta) y_ 2 } 
\cos \Big( \frac{\pi k}{2\ell} (y_2 + \ell) \Big) , \ j\in \bbZ,
\ k \in \bbN_0   
\end{eqnarray*}
The dispersion curves $\mu =  \mu_{jk}^0(\eta)$, $\eta \in [-\pi, \pi]$,
form the grid in Fig.\,\ref{fig4},a), where the disposition of its knots is
drawn in the case 
\begin{eqnarray*}
\frac16 < \ell < \frac14 . 
\end{eqnarray*}

We will consider the behavior of the dispersion curves of  the above described 
singular perturbation problem and show that the knots, marked with $\circ$
in Fig.\,\ref{fig4},b), disintegrate and give rise to small spectral gaps of width
$O(\varepsilon^2)$, which are presented by thick lines on the ordinate
axis.  

Asymptotic analysis as used here was initiated \cite{na453} and it has been
applied in different spectral problems with singular and regular
perturbation of boundaries, see \cite{BoPa,na530} and others.
However, in view of non-standard integro-differential boundary conditions
\eqref{38} on the small contour $\partial \theta^\varepsilon$, explicit 
formulas for the asymptotics of eigenvalues $\mu_{jk}^\varepsilon (\eta)$
would be required although error estimates would not be needed: the 
problem in $\omega \setminus \overline{\theta^\varepsilon}$ has a 
variational formulation, and the justification scheme follows word-by-word
those in \cite{BoPa,na453,na530}.

\subsection{Asymptotic ans\"atze} \label{sec4.2}
We consider the perturbation of the knot $(0 ,4 \pi^2)$ as in 
Fig.\,\ref{fig4},a). The eigenvalue $\mu_{\pm 10}(0) = 4 \pi^2$ of the problem 
\eqref{X2} has  multiplicity two and eigenfunctions are
\begin{eqnarray}
v_{\pm 10}(y) = e^{ \pm 2 \pi i y_1} .  \label{X5}
\end{eqnarray}
Following \cite{na453,na530} we define  the fast Floquet variable
\begin{eqnarray*}
\psi = \varepsilon^{-2} \eta     
\end{eqnarray*}
and introduce the asymptotic representation
\begin{eqnarray}
\mu^\varepsilon (\eta) = 4 \pi^2 + \varepsilon^2 \mu' (\psi) + \ldots \ . 
\label{X7}
\end{eqnarray}
The corresponding eigenfunctions are sought in the form
\begin{eqnarray*}
\ \ \ v^\varepsilon (\eta; y ) &=& v^0(\psi ; y ) + \varepsilon^2 v^2(\psi; 
u  ) + \varepsilon\chi_\theta(y) \big( w^0(\psi; \varepsilon^{-1} y  ) +
\varepsilon w^1 (\psi; \varepsilon^{-1} y  ) \big) + \ldots    
\end{eqnarray*}
where $\chi_\theta \in C^\infty( \overline \omega ) $ is a cut-off function
such that 
\begin{eqnarray*}
& & \mbox{ $\chi_\theta (y) = 1 $, if $|y_1| < 1/6$ and $|y_2| < \ell/3$,} 
\row
\mbox{  $\chi_\theta (y) = 0 $, if $|y_1| > 1/3$ or $|y_2| >  2\ell/3$},
\end{eqnarray*}
$v_\pm^p$ are  solutions of regular type in $\omega$ and $w^q$ are
 boundary layers written in the stretched coordinates $\zeta$,
see \eqref{X1}, in the aperture domain $\Upsilon = \bbR^2 \setminus
\overline{\theta^1}$.

The first regular term consists of the linear combination
\begin{eqnarray}
v^0(\psi; y  ) = A(\psi) e^{+ 2 \pi i y} + B(\psi) e^{- 2 \pi i y_1}
\label{X9}
\end{eqnarray}
where the  coefficients are to be determined. The second term $v_\pm^2$
satisfies the equations
\begin{eqnarray}
& & - \Delta_ y v^2 (\psi; y  ) - 4 \pi^2 v^2 (\psi; y ) = 
\mu' (\psi) v^0(\psi; y ) + f(\psi; y ), \  \ y \in \omega, 
\\
& &  \frac{\partial v^2 }{\partial y_2} (\psi; y_1, + \ell  ) 
- \frac{\partial v^2 }{\partial y_2} (\psi; y_1, - \ell  )  = 0, 
\ \ |y_a| < \frac12, \label{X10}
\end{eqnarray}
while, due to the relation $e^{i \eta} = e^{i \varepsilon^2 \psi}
= 1 + i \varepsilon^2 \psi + O(\varepsilon^4 \psi^2)$, the quasi-periodicity
conditions become
\begin{eqnarray}
& &   v^2  \Big(\psi; \frac12, y_2  \Big)   =
v^2  \Big(\psi; - \frac12, y_2  \Big) + i \psi v^0
\Big(\psi; - \frac12, y_2  \Big) ,  \ \ \ |y_2| < \ell
\row
\frac{\partial v^2 }{\partial y_1}  \Big(\psi; \frac12, y_2  \Big)   =
\frac{\partial v^2 }{\partial y_1}   \Big(\psi; - \frac12, y_2  \Big)
+ i \psi \frac{\partial v^0 }{\partial y_1} 
\Big(\psi; - \frac12, y_2  \Big), \ \ |y_2| < \ell. 
\label{X12}
\end{eqnarray}
The right-hand side $f$ will be determined in the next section 
after the examination of boundary layers. Then, the two compatibility
conditions in problem \eqref{X10}--\eqref{X12} will give us the 
correction term in \eqref{X7} as well as coefficients of the linear
combination \eqref{X9}.

\subsection{Boundary layer} \label{sec4.3}
Stretching the coordinates $y \mapsto \zeta = \varepsilon^{-1} y$ and
setting $\varepsilon=0$ turn the equations \eqref{25} with $H=1$ and
\eqref{38} with $\theta = \theta^\varepsilon$ into
\begin{eqnarray}
& & - \Delta_\zeta w(\zeta) = 0, \  \zeta \in \bbR^2 \setminus 
\overline{\theta^1} , \label{Y1}  \\
& & w(\zeta) = a \in \bbR, \  \zeta \in \partial \theta^1, 
\ \ \ \ \int\limits_{ \partial \theta^1} \partial_\nu w (\zeta) ds_\zeta
= 0 . \label{Y2}
\end{eqnarray}
Indeed, we have 
\begin{eqnarray*}
\Delta_y + \mu^\varepsilon(\eta) = \varepsilon^{-2}
\big( \Delta_\zeta + \varepsilon^2 \mu^\varepsilon(\eta) \big),
\ \ \ \mu^\varepsilon(\eta) {\rm mes}_2 \theta^\varepsilon
= O(\varepsilon^2)
\end{eqnarray*}
and, therefore, in the limit $\varepsilon \to +0$ the Helmholtz 
operator becomes the Laplacian and the right hand side of 
the integral condition \eqref{38} vanishes. 

We write the Taylor formula
\begin{eqnarray*}
& & v^0(\psi;y) = v^0(\psi;0) + y \cdot \nabla_y v^0(\psi;0) + O(|y|^2) 
\roweq 
v^2(\psi;0) + \varepsilon\zeta \cdot \nabla y v^0(\psi;0) + 
O(\varepsilon^2 |\zeta|^2 ) . 
\end{eqnarray*}
Since the constant $a$ is arbitrary, the first term $v^0(0)$
does not leave a discrepancy in the problem \eqref{Y1}, \eqref{Y2}.
Furthermore,
\begin{eqnarray*}
\int\limits_{ \partial \theta^1} \partial_\nu \zeta_j ds_\zeta
= \int\limits_{ \partial \theta^1} \nu_j(\zeta) ds_\zeta
= 0, \ \ \ j=1,2, 
\end{eqnarray*}
and therefore we need to find a harmonic function $\bfw$ satisfying
\begin{eqnarray}
\bfw_j(\zeta) = a_j - \zeta_ j, \ \ \zeta \in \partial \theta^1. 
\label{Y4}
\end{eqnarray}
We recall the definition of the polarization matrix $P(\theta^1)$
\cite[Appendix\, G]{PoSe}, which is a $(2 \times 2)$-matrix (a 
tensor of rank 2) and is composed from the coefficients in the 
decomposition of the decaying solution of the problem \eqref{Y1}, \eqref{Y4},
\begin{eqnarray}
\bfw_j(\zeta_j) = \frac{1}{2 \pi} \sum_{k=1}^2 P_{jk} (\theta^1) 
\frac{\zeta_k}{|\zeta|^2} + \widetilde \bfw_j (\zeta) ,  \ \ 
\widetilde \bfw_j (\zeta) = O\Big( \frac{1}{|\zeta|^2 } \Big) 
\label{Y5}
\end{eqnarray}
The matrix $P(\theta^1)$  is symmetric and positive definite, 
see \cite{PoSe}, and the last condition  in \eqref{Y2} is 
fulfilled due to Green's formula
\begin{eqnarray*}
\int\limits_{\partial \theta^1} \partial_\nu \bfw_j (\zeta) ds_\zeta
= - \lim_{R \to + \infty} \int\limits_{\partial \bbB_ R}
\frac{\partial \bfw^j}{\partial |\zeta|} (\zeta) ds_\zeta = 0 
\end{eqnarray*}
and the decay rate $O(|\zeta|^{-2})$ of the gradient $\nabla_\zeta
\bfw_j(\zeta) $, see \eqref{Y5}. Thus, this leads us to set
\begin{eqnarray}
w^1(\psi;\zeta ) = \sum_{k=1}^2 \bfw_j (\zeta) \frac{\partial v^0}{
\partial y_j} ( \psi; 0) . \label{Y6}
\end{eqnarray}

According to \cite{na61}, \cite[Ch.\,9,10]{MaNaPl} we also need to 
examine the next boundary layer term $w^2(\eta; \zeta)$ which takes into 
account the second-order term
\begin{eqnarray*}
Q(\psi;y ) = \frac12 \sum_{p,q=1}^2 y_p y_q \frac{\partial^2 v^0}{\partial
y_p \partial y_q} (\psi; 0) 
\end{eqnarray*}
in the Taylor formula for $v^0$. Noting that $|\theta^\varepsilon| = 
\varepsilon^2 |\theta^1|$ we get 
\begin{eqnarray*}
\ \ \ 
\int\limits_{\partial \theta^1} \partial_\nu Q(\psi;\zeta ) ds_\zeta 
= \int\limits_{\theta^1} \Delta_\zeta  Q(\psi;\zeta ) d \zeta 
= |\theta^1 | \Delta_y v^0 (\psi;0 )
= - 4 \pi^2 |\theta^1 | v^0(\psi;0 ).   
\end{eqnarray*}
Consequently, the functions $w^2$
must satisfy the Laplace equation as well as the conditions 
\begin{eqnarray}
& & w^2(\psi; \zeta ) = a^2 (\psi) - Q(\psi;\zeta  ) , 
\row
\int\limits_{\partial \theta^1} \partial_\nu w^2 (\psi;\zeta ; 
\psi) ds_\zeta = 4 \pi^2 | \theta^1 | v^0 (\psi;0,\psi) 
=: T^0(\psi)   \label{Y8}
\end{eqnarray}
and thus  it can be  decomposed as 
\begin{eqnarray*}
w^2(\psi;\zeta  ) = T^0 (\psi) \frac{1}{2 \pi} \ln \frac{1}{|\zeta|}
+ T^1(\psi ; \ln \varepsilon ) + \widetilde w^2 (\psi;\zeta  ),
\ \ \ 
\widetilde w^2 (\psi;\zeta  ) = O\Big( \frac1{|\zeta|} \Big) .
\end{eqnarray*}

We mention that the harmonic functions $w_\pm^2$ grow at infinity
and thus it is convenient to choose the constant
\begin{eqnarray*}
T^1(\psi ; \ln \varepsilon) = T^0 (\psi) \frac{1}{2 \pi}  \ln 
\frac1\varepsilon ,
\end{eqnarray*}
which allows us to write 
\begin{eqnarray}
\varepsilon^2 w^2 (\psi; \zeta  ) = \varepsilon^2 T^0(\psi) 
\frac{1}{2 \pi}  \ln \frac{1 }{|y|} + O\Big( \frac{\varepsilon^3}{|y|} \Big).
\label{Y9}
\end{eqnarray}
Finally, we take into account \eqref{Y6}, \eqref{Y5} and obtain
\begin{eqnarray}
\varepsilon w^1(\psi;\zeta ) = \frac{\varepsilon^2}{2 \pi}
\sum_{j,k=1}^2  P_{jk} (\theta^1) \frac{y_k}{|y_k|} 
\frac{\partial v^0}{\partial y_k} (\psi;0 ) + 
O\Big( \frac{\varepsilon^3}{|y|} \Big).  \label{Y10}
\end{eqnarray}
Now we insert the ansatz \eqref{Y8} into the equation \eqref{25},
separate terms of order $\varepsilon^2$, and taking into account 
\eqref{Y9}, \eqref{Y10}, complete the composition of the equation 
\eqref{X10} by setting
\begin{eqnarray*}
& & f(\psi; y ) =  \big( \big[ \Delta_y, \chi_\theta(y) \big] + \mu^0(\psi) 
\chi_\theta(y) \big) \cT(\psi;y  ), \ \ \mbox{where} \row
\cT(\psi; y ) = T^0(\psi) \frac{1}{2\pi } \ln \frac{1}{|y|} + 
\sum_{j,k=1}^2  P_{jk} (\theta^1) \frac{y_k}{|y_k|} 
\frac{\partial v^0}{\partial y_k} (\psi;0 ). 
\end{eqnarray*}

\subsection{Algebraic system for $\mu'(\psi)$}  \label{sec4.4}
Since $4 \pi^2$ is an eigenvalue of multiplicity 2, the problem 
\eqref{X10}--\eqref{X12}  must be associated with two compatibility
conditions, which can be derived by inserting a possible solution
$v_\pm^2$ and the eigenfunctions \eqref{X5} into Green's formula
\begin{eqnarray}
& & I_\pm^1(\psi) + I_\pm^2(\psi)  := \mu_\pm'(\psi) 
\int\limits_\omega \overline{v_{\pm10}(y)} v^0(\psi;y ) dy
+ \int\limits_\omega \overline{v_{\pm10}(y)} f_\pm (\psi;y ) dy
\roweq
\int\limits_{-\ell}^\ell \Big( v^2(\psi;y  ) 
\overline{\frac{\partial v_ {\pm 10}}{\partial y_ 1}(y)} - 
\overline{ v_ {\pm 10}(y) } \frac{\partial v^2}{\partial y_ 1}
(\psi;y ) \Big) \bigg|_ {y_1= -1/2}^{y_1= 1/2} dy_2  =: I_\pm^3(\psi) .
\label{T1}
\end{eqnarray}
Let us compute the integrals $I_\pm^p(\psi)$. Owing to \eqref{X5} 
and \eqref{X9}, we readily obtain
\begin{eqnarray*}
I_+^1(\psi)  =2 \ell A (\psi) \mu'(\psi) , \ \ 
I_-^1(\psi)  =2 \ell B (\psi) \mu'(\psi) .   
\end{eqnarray*}
Using Green's formula in the perforated rectangle $\omega \setminus 
\bbB_\delta$ yields
\begin{eqnarray*}
& & I_\pm^2(\psi) = - \lim\limits_{\delta \to +0} 
\int\limits_{\omega \setminus \bbB_\delta} \overline{v_{\pm 10} (y) } 
(\Delta_y + 4 \pi^2) \big( \chi_\delta(y) \cT(\psi;y ) \big) dy
\roweq
\lim\limits_{\delta \to +0}  \int\limits_{\omega \setminus \bbB_\delta} 
\Big( \cT(\psi; y ) \overline{\partial_{|y|} v_{\pm10}(y) }
-  \overline{ v_{\pm10}(y) } \partial_{|y|}  \cT(\psi;y ) \Big) ds_y
\roweq
4 \pi^2 \big( |\theta'| ( A(\psi) + B(\psi) \big) \pm
P_{11}(\theta')  \big( ( A(\psi) - B(\psi) \big) . 
\end{eqnarray*}
Recalling the relations \eqref{X12} gives us 
\begin{eqnarray}
& & I_\pm^3(\psi) = 4 \pi \ell \psi \Big( 
\big( A(\psi) + B(\psi) \big) +   \big( ( A(\psi) - B(\psi) \big) \Big)
\nonumber \\
\Rightarrow & &  I_+^3(\psi) = 8 \pi \ell \psi  A(\psi), \ \ \ 
I_-^3(\psi) = - 8 \pi \ell \psi  A(\psi). \label{T4}
\end{eqnarray}
Combining formulas \eqref{T1}--\eqref{T4} we arrive at the system
of linear algebraic equations
\begin{eqnarray*}
& & \big( \cP_+(\theta^1) - 8 \pi \ell \psi \big) A(\psi)  + 
\cP_-(\theta^1) B(\psi) = 
- 2 \ell A(\psi) \mu'(\psi) , 
\row
 \cP_-(\theta^1)  A(\psi)  + \big( 
\cP_+(\theta^1) +  8 \pi \ell \psi \big) B(\psi)  = 
- 2 \ell B(\psi) \mu'(\psi) , 
\end{eqnarray*}
where 
\begin{eqnarray*}
\cP_\pm (\theta^1) = 4 \pi^2 \big( |\theta^1| \pm P_{11}(\theta^1) \big) .
\end{eqnarray*}
Thus,
\begin{eqnarray*}
\mu_\pm'(\psi) = - \frac{1}{2 \ell} \Big( \cP_\pm(\theta^1)  \pm
\sqrt{\cP_\pm (\theta^1)^2 + 64 \pi^2 \ell^2 \psi^2 } \Big) .
\end{eqnarray*}

\subsection{Conclusions on spectral gaps} \label{sec4.5}
\begin{figure}
\begin{center}
\includegraphics[width=10cm]{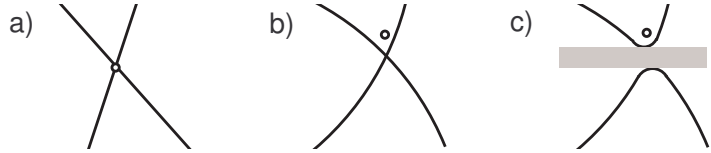}
\end{center}
\caption{a), b) Knot of the dispersion curves,  c) disintegrating knot.}
\label{fig5}
\end{figure}

If $\cP_-(\theta^1) =0$, the graphs of the functions $\psi \mapsto 
\mu_\pm'(\psi)$ include two crossing straight lines (Fig.\,\ref{fig5},a),
but in the case 
\begin{eqnarray}
\cP_-(\theta^1) \not= 0   \label{T7}
\end{eqnarray} 
the graphs turn into two disjoint parabolas (Fig.\,\ref{fig5},c)). 
Thus, in the same way as in the papers cited above, we conclude that the knot 
$(4 \pi^2,0)$ in Fig.\,\ref{fig4},b) disintegrates as in Fig.\,\ref{fig4},c), 
and causes a  spectral gap of width $\ell^{-1} \varepsilon^2 
\cP_-(\theta^1)) + O(\varepsilon^3)$ with the center at the point 
$4 \pi^2 - (2 \ell)^{-1} \varepsilon^2 \cP_+(\theta^1) + O(\varepsilon^3)$.

To give an example where \eqref{T7} holds, we observe that for a horizontal 
crack  $\overline{\theta^1 } = \{ \zeta \, : \, \zeta_2 = 0, \ 
|\zeta_1 | < {\sf L} \}$ we have
\begin{eqnarray*}
|\theta' |= 0 \ \ \mbox{and} \ \ P_{11} (\theta') = \pi {\sf L}^2 > 0 , 
\end{eqnarray*}
see \cite[Appendix G]{PoSe}. 

An example of a smooth domain with \eqref{T7} is given by observing
that for a thin ellipse $\theta_\delta^1$ with axes $1 $ and $\delta$
we have (see again \cite[Appendix G]{PoSe})
\begin{eqnarray*}
|\theta_\delta^1  |=  O(\delta)  \ \ \mbox{and} \ \ 
P_{11} (\theta_\delta^1 ) = \pi^2 ( 1 + O(\delta^2) )/4 > 0 .
\end{eqnarray*}
.

The equality $\cP_-(\theta) = 0$ does not yet imply the non-existence of
a gap, but this depends on  higher  order asymptotic terms. 

Similar calculations can be applied to study the existence of gaps
near all knots marked with $\circ$ in Fig.\,\ref{fig4},b). On the other hand, the 
knots marked with $\bullet$, which are obtained as the intersection of 
either two ascending or two descending curves, never get disjoined
(see an explanation in, e.g., \cite{na453}).


\section*{Conflict of interest}

The authors declare that they have no conflict of interest.

\end{document}